%

\documentstyle{amsppt}
\UseAMSsymbols
\catcode`\@=11

\def\leftrightarrowfill{$\m@th\mathord\leftarrow\mkern-6mu%
  \cleaders\hbox{$\mkern-2mu\mathord-\mkern-2mu$}\hfill
  \mkern-6mu\mathord\rightarrow$}

\atdef@?#1?#2?{\ampersand@\setbox\z@\hbox{$\ssize
 \;\;{#1}\;$}\setbox\@ne\hbox{$\ssize\;\;{#2}\;$}\setbox\tw@
 \hbox{$#2$}\ifCD@
 \global\bigaw@\minCDaw@\else\global\bigaw@\minaw@\fi
 \ifdim\wd\z@>\bigaw@\global\bigaw@\wd\z@\fi
 \ifdim\wd\@ne>\bigaw@\global\bigaw@\wd\@ne\fi
 \ifCD@\hskip.5em\fi
 \ifdim\wd\tw@>\z@
 \mathrel{\mathop{\hbox to\bigaw@{\leftrightarrowfill}}\limits^{#1}_{#2}}\else
 \mathrel{\mathop{\hbox to\bigaw@{\leftrightarrowfill}}\limits^{#1}}\fi
 \ifCD@\hskip.5em\fi\ampersand@}
\atdef@-#1-#2-{\ampersand@\setbox\z@\hbox{$\ssize  
 \;\;{#1}\;$}\setbox\@ne\hbox{$\ssize\;\;{#2}\;$}\setbox\tw@
 \hbox{$#2$}\ifCD@
 \global\bigaw@\minCDaw@\else\global\bigaw@\minaw@\fi
 \ifdim\wd\z@>\bigaw@\global\bigaw@\wd\z@\fi
 \ifdim\wd\@ne>\bigaw@\global\bigaw@\wd\@ne\fi
 \ifCD@\hskip.5em\fi
 \ifdim\wd\tw@>\z@
 \mathrel{\mathop{\raise.5ex\hbox to\bigaw@{\hrulefill}}\limits^{#1}_{#2}}
\else\mathrel{\mathop{\raise.5ex\hbox to\bigaw@{\hrulefill}}\limits^{#1}}\fi
 \ifCD@\hskip.5em\fi\ampersand@}

\def\hookrightarrowfill{$\m@th\mathord\lhook\mkern-3.1mu%
  \cleaders\hbox{$\mkern-2mu\mathord-\mkern-2mu$}\hfill
  \mkern-6mu\mathord\rightarrow$}
\atdef@(#1(#2({\ampersand@\setbox\z@\hbox{$\ssize
 \;\;{#1}\;$}\setbox\@ne\hbox{$\ssize\;\;{#2}\;$}\setbox\tw@
 \hbox{$#2$}\ifCD@
 \global\bigaw@\minCDaw@\else\global\bigaw@\minaw@\fi
 \ifdim\wd\z@>\bigaw@\global\bigaw@\wd\z@\fi
 \ifdim\wd\@ne>\bigaw@\global\bigaw@\wd\@ne\fi
 \ifCD@\hskip.5em\fi
 \ifdim\wd\tw@>\z@
 \mathrel{\mathop{\hbox to\bigaw@{\hookrightarrowfill}}\limits^{#1}_{#2}}\else
 \mathrel{\mathop{\hbox to\bigaw@{\hookrightarrowfill}}\limits^{#1}}\fi
 \ifCD@\hskip.5em\fi\ampersand@}

\parindent20\p@
\advance\captionwidth@1in

\font\chapheadfont@=cmr10 scaled\magstep2
\font\chapheadmathfont@=cmmi10 scaled\magstep2
\font\chapheadmsa=msam10 scaled\magstep2
\font\chapheadmsb=msbm10 scaled\magstep2
\font\chapheadscriptmathfont=cmmi7 scaled\magstep2
\define\chapheadfont{\chapheadfont@\textfont1=\chapheadmathfont@%
        \scriptfont1=\chapheadscriptmathfont%
        \textfont\msafam=\chapheadmsa\textfont\msbfam=\chapheadmsb}
\outer\def\chapheading{\newpage
  \begingroup\raggedcenter@\interlinepenalty\@M \let\\\linebreak
  \chapheadfont\noindent\ignorespaces}

\newif\ifdatver
\define\datver#1{\ifdatver\global\def\datverp{}\else\datvertrue%
  \global\def\datverp{\par\boxed{\boxed{\text{Version: #1; Run: \today}}}}\fi}

\newif\ifbigdoc

\let\section\relax

\def\nosection{\newcodes@\endlinechar=10 \sect@}
{\lccode`\!=`\\
\lowercase{\gdef\sect@#1^^J{\sect@@#1!section\sect@@@}%
\gdef\sect@@#1!section{\futurelet\next\sect@@@}%
\gdef\sect@@@#1\sect@@@{\ifx\next\sect@@@\let
\next=\sect@\else\def\next{\oldcodes@\endlinechar=`\^^M\relax}%
 \fi\next}}}
\catcode`\@=\active
\define\protag#1 #2{{\hbox{\rm\ignorespaces#1\unskip}}#2}
\define\theprotag#1 #2{#2 {\rm\ignorespaces#1\unskip}}
\define\exertag #1 #2{\demo{\hbox{\rm(\ignorespaces#1\unskip)} #2}}

\define\figtagg#1{Figure $#1$}

\NoRunningHeads
\CenteredTagsOnSplits
\NoBlackBoxes

\def\today{\ifcase\month\or
 January\or February\or March\or April\or May\or June\or
 July\or August\or September\or October\or November\or December\fi
 \ \space\number\day, \number\year}

\define\Dom{\operatorname{Dom}}

\define\OS {\negthinspace\smallsetminus\negthinspace 0}

\define\dCc#1{\dot{\Cal{C}}^{#1}_c}

\define\nov#1{{\frac{n}{#1}}}

\define\sgn{\operatorname{sgn}}

\define\({\bigl(}
\define\){\bigr)}
\magnification=\magstep1
\hsize 6.0truein
\vsize 8.2truein
\hcorrection{.2truein}
\vcorrection{.2truein}
\loadbold       

\define\sg{\Sigma^*_n}
\redefine\sgn{\Sigma^*_{n-1}}
\define\so{\Sigma^*_0}

\define\ORD{\operatorname{ORD}}
\define\Rng{\operatorname{Rng}}
\define\card{\operatorname{card}}
\define\Card{\operatorname{Card}}
\define\dom{\operatorname{dom}}
\define\cof{\operatorname{cof}}

\define\Lim{\operatorname{Lim}}

\define\even{{\text{even}}}
\define\odd{{\text{odd}}}
\redefine\phi{\varphi}
\document
\baselineskip=15pt

\font\bigtenrm=cmr12 scaled\magstep2
\centerline
{\bigtenrm {A Simpler Proof of Jensen's Coding Theorem}}
\vskip20pt

\font\bigtenrm=cmr10 scaled\magstep2
\centerline
{\bigtenrm{Sy D. Friedman}\footnote"*"{Research supported by NSF Grant 
\# 9205530.}}

\centerline
{\bigtenrm {M.I.T.}}

\vskip20pt

\comment
lasteqno 1@0
\endcomment

Beller-Jensen-Welch [82] provides a proof of Jensen's remarkable Coding
Theorem, which demonstrates that the universe can be included in $L[R]$ for
some real $R,$  via class forcing. The purpose of this article is to present a
simpler proof of Jensen's theorem, obtained by implementing some changes first
developed for the theory of strong coding (Friedman [87]).

The basic idea is to first choose $A\subseteq ORD$ so that $V=L[A]$ and then
generically add sets $G_\alpha\subseteq\alpha^+,\alpha$ $O$ or an infinite
cardinal ($O^+$ denotes $\omega)$  so that $G_\alpha$  codes both
$G_{\alpha^{+}}$ and $A\cap\alpha^+.$  Also for limit cardinals $\alpha,
G_\alpha$ is coded by $\langle G_{\bar\alpha }|\bar\alpha<\alpha\rangle.$
Thus there are two ``building blocks'' for the forcing, the successor coding
and the limit coding.  We modify the successor coding so as to eliminate
Jensen's use of ``generic codes'' (this improves an earlier modification of 
this type, due to Welch and Donder). And we thin out the limit coding so as to
eliminate the technical problems causing Jensen's split into cases according
to whether or not $O^{\#}$ exists.

\proclaim{Theorem} (Jensen) \  There is a class forcing ${\Cal{P}}$ such that
if $G$  is ${\Cal{P}}$-generic over $V$  then $V[G]\models ZFC+V=L[R],
R\subseteq\omega.$  If $V\models GCH$ then ${\Cal{P}}$  preserves cardinals.
\endproclaim

It is not difficult to class-generically extend $V$  to  make $GCH$ true. And
any ``reshaped'' 
subset of $\omega_1$  can be coded by a real via a $CCC$  forcing. 
(See Section One below for a definition of ``reshaped''.)
So it
suffices to prove that $V$  can be coded by a ``reshaped'' 
subset of $\omega_1,$
preserving cardinals, assuming the $GCH.$  As a first step, force $A\subseteq
ORD$  such that for each infinite cardinal $\alpha,$  $L_\alpha [A]=H_\alpha=$
all sets of hereditary cardinality less than $\alpha.$ 

\vskip10pt

\flushpar
{\bf Section One \ The Successor Coding $R^s.$ }

Fix an infinite cardinal $\alpha.$  $S_\alpha$ is defined to be a certain
collection of ``strings'' $s:\ [\alpha,|s|)\longrightarrow 2,\alpha \le$ 
$|s|<\alpha^+.$
For $s$  to belong to $S_\alpha$  we require that $s$  is ``reshaped''. This
means that for $\eta\le|s|,$ $L[A\cap\alpha, s\restriction\eta
]\models\card(\eta)\le\alpha.$  The reshaping of $s$ allows us to code $s$  by
a subset of $\alpha,$  in the manner which we now describe.

For $s\in S_\alpha$  define structures
${\Cal{A}}^0_s=L_{\mu^{0}_{s}}[A\cap\alpha,s^*],$
${\Cal{A}}_s=L_{\mu_s}[A\cap\alpha,s^*]$  as follows (where
$s^*=\{\mu_{s\restriction\eta }|s(\eta)=1\}):$ If $|s|=\alpha$  then
$\mu^0_s=\alpha.$  For $|s|>\alpha,$ $\mu^0_s=\bigcup\{\mu_{s\restriction\eta
}|\eta<|s|\}$  and in general $\mu_s=$ least p\.r\. closed ordinal $\mu$
greater than $\mu^0_s$ such that $L_\mu
[A\cap\alpha,s^*]\models\card(s)\le\alpha.$  These ordinals are well-defined
due to the reshaping of $s.$  

For $s\in S_\alpha$ we write $\alpha(s)=\alpha.$  Note that if $|s|=\alpha(s)$
then $s=\emptyset;$ in this case we think of $s$  as ``labelled'' with the
ordinal $\alpha(s),$  so that there are distinct $s_\alpha\in
S_\alpha,\alpha(s_\alpha)=\alpha.$  

For later use we also define structures $\widehat{\Cal{A}}_s$ and
${\Cal{A}}'_s$  for $s\in S_\alpha$  as follows: let $\hat\mu_s=$ largest
p\.r\. closed $\mu$  such that $\mu=\mu^0_s$ or $L_\mu
[A\cap\alpha,s^*]\models|s|$  is a cardinal greater than $\alpha.$  Then 
$\widehat{\Cal{A}}_s=L_{\hat\mu_s}[A\cap\alpha,s^*].$  The ordinal $\mu'_s$
and structure ${\Cal{A}}'_s$  are defined in the same way, except we replace
p\.r\. closure of $\mu$  by the weaker condition $\omega\cdot\mu=\mu.$ 

For $s\in S_{\alpha^+}$ write $\bar s<s$ to mean that $\pi(\bar s)=s$ where
$\pi:\ \overline{\Cal{A}}\longrightarrow {\Cal{A}}_s$ is an elementary
embedding with some critical point $\alpha(\bar s)<\alpha^+$  and where
$\pi(\alpha(\bar s))=\alpha^+.$  Then $\pi=\pi_{\bar ss}$ is unique. Let $\bar
s\le s$ denote $\bar s<s$ or $\bar s=s.$  We have the following facts:

(a) \  $\{\alpha(\bar s)|\bar s<s\}$ is $CUB$ in $\alpha^+.$

(b) \  If $\bar t$ is a proper initial segment of $\bar s$  then $\bar
t<\pi_{\bar ss}(\bar t)=t$  and $\pi_{\bar tt}=\pi_{\bar ss}\restriction
{\Cal{A}}_{\bar t}.$

(c) \ ${\Cal{A}}_s=\bigcup\{\Rng(\pi_{\bar ss})|\bar s<s\}.$

Now for $s\in S_{\alpha^+}$ let $b_s=\{\bar s|\bar s<s\}.$  We use the strings
$\bar s*i$ with $\bar s
<s\restriction\eta,$ $i=0$ or $1,$  to code $s(\eta).$  A
condition in the successor coding $R^s$ is a pair $(u,\bar u)$ where:

1) \   $u\in S_\alpha$ 

2) \  $\bar u\subseteq\{b_{s\restriction\eta}|s(\eta)=0\},$ card $(\bar
u)\le\alpha$ in ${\Cal{A}}_s.$ 

\flushpar
To define extension of conditions, we need a couple of preliminary
definitions. We say that $\bar u$ {\it restrains} $\bar s* 1$ if $\bar s\in
b$ for some $b\in\bar u$ and $\bar s$ {\it lies on} $u$ if
$u(\alpha(\bar s))=1$ and $u(\langle\alpha(\bar s),\eta\rangle)=\bar
s(\eta)$
for $\eta\in$ Dom$(\bar s).$
  Also let $\langle Z_\gamma|\gamma<\alpha^+\rangle$
be an $L_{\alpha^{+}}$-definable partition of the odd ordinals less than
$\alpha^+$ into $\alpha^+$ disjoint pieces of size $\alpha^+.$  We use the
$Z_\gamma$'s to code $A\cap\alpha^+$ into $G_\alpha.$  For $u\in S_\alpha,$
$u^{\text{even}}$ $(\delta)=u(2\delta),$  $u^{\text{odd}}$
 $(\delta)=u(2\delta+1).$ 

Extension of conditions for $R^s$ is defined by:\ $(u_0,\bar u_0)\le(u_1,\bar
u_1)$ iff $u_0$ extends $u_1;\bar u_0\supseteq\bar u_1;$  $\bar u_1$ restrains
$\bar s*1,$ $\bar s*1$ lies on $u_0^{\text{even}}\longrightarrow\bar
s*1$ lies on $u_1^{\text{even}};$ $\gamma<|u_1|,\gamma\notin A,$ $\delta\in
Z_\gamma,$ $u_0^{\text{odd}}(\delta)=1\longrightarrow
u_1^{\text{odd}}(\delta)=1.$  Note that $R^s\in{\Cal{A}}_s.$

\proclaim{Lemma 1.1} Suppose $G$  is $R^s$-generic over ${\Cal{A}}_s$  and let
$G_\alpha=\bigcup\{u|(u,\bar u)\in G$ for some $\bar u\}.$  Then
$G,A\cap\alpha^+,s$ belong to $L_{\mu_{s}}[G_\alpha ].$ 
\endproclaim

\demo{Proof} We can write $(u,\bar u)\in G$ iff $u\subseteq G_\alpha$  and
$\bar s\in b\in\bar u,$ $\bar s*1$ lies on $\ G_\alpha^{\text{even}}\longrightarrow \bar 
s*1$ lies on $\ u^{\text{even}}$ and $\gamma<|u|,$
$\gamma\notin A,\delta\in 
Z_\gamma,$  $G_\alpha^{\odd}(\delta)=1\longrightarrow u^{\odd}(\delta)=1.$  So
$G\in L_{\mu_{s}}[A\cap\alpha^+,G_\alpha,s ].$  And $\gamma\in A\cap\alpha^+$
iff $G^{\odd}_\alpha(\delta)=1$  for unboundedly many $\delta\in Z_\gamma,$
so $G,A\cap\alpha^+\in L_{\mu_{s}}[G_\alpha,s].$  Finally note that for any
$\eta<|s|,\bar s$ lies on $G^{\even}_\alpha$ for unboundedly many $\bar
s<s\restriction\eta$ by a density argument using the fact that for
$\eta<|s|,(u,\bar u)\in R^s,$ $b_{s\restriction\eta}$ is almost disjoint
from $\{u|u$ extends some $\bar s*1$ restrained by $\bar u\}$.  
So $s(\eta)=1$ iff
$\bar s*1$ lies on $G_\alpha^{\even}$ for unboundedly many $\bar
s<s\restriction\eta.$  Thus $s\restriction\eta$ can be recovered by
induction on $\eta\le|s|,$ inside $L_{\mu_{s}}[G_\alpha ].$ 
\hfill{$\dashv$ }
\enddemo

\proclaim{Lemma 1.2} $R^{<s}=\bigcup\{R^t|t\subseteq s, t\ne s\}$ has the
$\alpha^{++}$-CC in $\widehat{\Cal{A}}_s.$ 
\endproclaim

\demo{Proof} If $\hat\mu_s=\mu^0_s$ then this is vacuous. Otherwise we need
only observe that $R^{<s}\in\widehat{\Cal{A}}_s$ and $(u_0,\bar u_0),$
$(u_1,\bar u_1)$ incompatible $\longrightarrow u_0\ne u_1$ and $S_\alpha$ has
cardinality $\alpha^+$ in $\widehat{\Cal{A}}_s.$  \hfill{$\dashv$ }
\enddemo

\proclaim{Lemma 1.3}  $R^s$ is $\le\alpha$-distributive in ${\Cal{A}}_s.$ 
\endproclaim

\demo{Proof} Suppose $(u_0,\bar u_0)\in R^s$ and $\langle D_i|i<\alpha\rangle$
are predense on $R^s,$ $\langle D_i|i<\alpha\rangle\in {\Cal{A}}_s.$   By 
induction we define conditions $(u_i,\bar u_i)$ and elementary submodels $M_i$
of ${\Cal{A}}_s$ with $(u_i,\bar u_i)\in M_{i+1},$  for $i\le\alpha.$  Choose
$M_0$  to contain $\alpha$  as a subset and to contain $\langle
D_i|i<\alpha\rangle,s,A\cap\alpha^+$  as elements. Having defined $(u_i,\bar
u_i)$ and $M_i,$  choose $M_{i+1}$  to contain $M_i$  as a subset and
$(u_i,\bar u_i)$  as an element.  Choose $(u_{i+1},\bar u_{i+1})$  to extend
$(u_i,\bar u_i),$  meet $D_i,$  guarantee that if $s(\eta)=1,$ $\eta\in M_i$
then $\bar s*1$ lies on $u^{\even}_{i+1}-u^{\even}_i$ for some $\bar
s<s\restriction\eta,$  guarantee that if $\gamma\in
A\cap(M_i\cap\alpha^+)$  then $u^{\odd}_{i+1}(\delta)=1$  for some
$\delta\notin\dom u_i^{\text{odd}},$ $\delta \in Z_\gamma,$  and finally
choose $\bar u_{i+1}$  to 
contain all $b_{s\restriction\eta}$ with $s(\eta)=0,$ $\eta\in M_i.$  The
last requirement can be imposed because the facts that $|s|$  has cardinality
$\le\alpha^+$ in ${\Cal{A}}_s,$  $H_{\alpha^{+}}\subseteq {\Cal{A}}_s$  imply
that any subset of $|s|$  of cardinality $\le\alpha$  belongs to
${\Cal{A}}_s.$ 

For $\lambda\le\alpha$ limit, $M_\lambda=\bigcup\{M_i|i<\lambda\}$ and
$u_\lambda=\bigcup\{u_i|i<\lambda\},$ $\bar u_\lambda=\bigcup\{\bar
u_i|i<\lambda\}.$  By construction, $u_\lambda$ codes
$A\cap(M_\lambda\cap\alpha^+)$  as well as $\bar s=s\circ\pi^{-1}$ where $\pi$
is the transitive collapse map for $M_\lambda.$  Thus the sequence of ordinals
$\langle M_i\cap\alpha^+|i<\lambda\rangle$ is cofinal in $|u_\lambda|$ and
belongs to $L[u_\lambda ],$ since the entire sequence
$\langle\overline{M}_i|i<\lambda\rangle$  can be recovered in $L[u_\lambda ],$
$\overline{M}_i=$ transitive collapse $(M_i).$  This shows that $u_\lambda$ is
reshaped, so $(u_\lambda,\bar u_\lambda)$ is a condition. Finally note that 
$(u_\alpha,\bar u_\alpha)$ is an extension of $(u_0,\bar u_0)$ meeting each of
the $D_i$'s. \hfill{$\dashv$ }

\enddemo

\proclaim{Corollary 1.4} $R^{<s}$ is $\le\alpha$-distributive in
$\widehat{\Cal{A}}_s.$  
\endproclaim

\demo{Proof} By Lemma 1.2 it suffices to prove $\le\alpha$-distributivity
in ${\Cal{A}}^0_s.$  This is easily proved by induction on $|s|,$ using Lemma
1.3 at successor stages. \hfill{$\dashv$ } 
\enddemo

\proclaim{Lemma 1.5} If $D\subseteq R^{<s},$ $D\in\widehat{\Cal{A}}_s$ is 
predense and $s\subseteq t\in S_{\alpha^{+}}$  then $D$  is predense on 
$R^t.$ 
\endproclaim

\demo{Proof} It suffices to show that if $D\subseteq R^s,$ $D\in{\Cal{A}}_s$
is predense, $s\subseteq t\in S_{\alpha^{+}}$ then $D$  is predense on $R^t;$
for then, as in the proof of Corollary 1.4, we can induct on $|s|$  and use
Lemma 1.2.

Suppose $D$ is predense on $R^s,D\in{\Cal{A}}_s$ and $(u,\bar u)$  belongs
to $R^t.$  We can extend $(u,\bar u)$  to guarantee  that for some $\bar
t<t,\bar u=\{b_{t\restriction(\eta+1)}|t(\eta)=0,\eta\in\Rng\pi_{\bar tt}\},$
$D,s\in\Rng(\pi_{\bar tt})$ and $|u|=\alpha(\bar t)+1,$ $u(\alpha(\bar
t))=0,$
$u\in\Cal{A}_{\bar t\restriction\alpha(\bar t)}$.
 Let $(u^*,\bar u^*)$ be the least extension of
$(u,\bar u\cap{\Cal{A}}_s)\in R^s$ meeting $D.$  We claim that $(u^*,\bar
u^*\cup\bar u)$ is an extension of $(u,\bar u),$  and this will prove the
lemma. Clearly $\gamma<|u|,$ $\delta\notin A,$ $\delta \in Z_\gamma,$
$u^{*\odd}(\delta)=1\longrightarrow u^{\odd}(\delta)=1,$ since $(u^*,\bar
u^*)$ extends $(u,\bar u\cap{\Cal{A}}_s).$  Suppose
$r<t\restriction\eta, t(\eta)=0$  where $\eta\in\Rng\pi_{\bar tt}$
and $r*1$ lies 
on $u^{*\even}.$  If $\eta<|s|$  then $r*1$ lies on $u^{\even},$ as desired,
since $(u^*,\bar u^*)$ extends $(u,\bar u\cap{\Cal{A}}_s).$  If
$\alpha(r)<\alpha(\bar t)$ then $|r|<\alpha(\bar t)$ so again $r*1$ lies
on $u^{\even}$ since $|u|>\alpha(\bar t)>|r*1|.$  If
$\alpha(r)=\alpha(\bar t)$  then 
 $r*1$ cannot lie on $u^{*\even},$ by choice of $u.$
  Finally if $\alpha(r)>\alpha(\bar t)$ then since
$\eta\ge|s|$ we have $\alpha(r)>|u^*|$ by leastness of $(u^*,\bar u^*).$  So
$r*1$ cannot lie on $u^{*\even}.$ \hfill{$\dashv$ }
\enddemo

\vskip20pt

\flushpar
{\bf Section Two \  Limit Coding.}

We begin with a rough indication of the forcing ${\Cal{P}}^u$  for coding
$u\in S_\alpha,\alpha$  an uncountable limit cardinal, into a subset of
$\alpha.$  ${\Cal{P}}^u\subseteq {\Cal{A}}_u$ consists of
${\Cal{P}}^{<u}=\bigcup\{{\Cal{P}}^{u\restriction\xi }|\xi<|u|\}$ together
with certain $p:\Card\cap\alpha\longrightarrow V$ such that
$p(\beta)=(p_\beta,\bar p_\beta)\in R^{p_{\beta^{+}}}$ for $\beta\in\dom(p).$
(We use Card to denote the class of infinite cardinals.) Also for uncountable
limit cardinals $\beta<\alpha$  we (inductively) require that
$p\restriction\beta\in {\Cal{P}}^{p_{\beta }}-{\Cal{P}}^{<p_{\beta }}.$  We
also insist that $p$  code $u$ in the following sense: For $\xi<|u|$  and
$\beta\in\Card\cap\alpha$ define $M^\xi_\beta=\Sigma_1$ Skolem hull of
$\beta\cup\{u\restriction\xi,A\cap\alpha\}$ in ${\Cal{A}}_{u\restriction\xi }$
and $b^\xi_\beta=M^\xi_\beta\cap\beta^+.$  Then code $u$  by: $u(\xi)=1$ iff
$p^{\odd}_{\beta^{+}}(b^\xi_{\beta^{+}})=1$ for sufficiently large
$\beta\in\Card\cap\alpha.$  Recall that the successor coding
$R^{p_{\beta^{++}}}$  makes use of odd ordinals (in the $Z_\gamma$'s) so
the successor and limit codings do not conflict.  For $p,q\in {\Cal{P}}^u$ we
write $p\le q$ iff $p(\beta)\le q(\beta)$ in $R^{p_{\beta^{+}}}$ for each
$\beta\in\Card\cap\alpha.$ 

To facilitate the proofs of extendibility and distributivity for 
${\Cal{P}}^u$  we thin out the forcing, in a number of ways. For this purpose
we need appropriate forms of $\square$ and $\diamond,$ in a relativized form.
Jensen  observed that his proofs of these principles for $L$ go through when
relativized to reshaped strings. Precisely:

\flushpar
{\bf Relativized $\bold{\square}$ } \ Let $S=\bigcup\{S_\alpha|\alpha$ an
infinite cardinal$\}.$  There exists  $\langle C_s|s\in S\rangle$ such that
$C_s\in{\Cal{A}}_s$ and:

(a) \  If $\alpha(s)<|s|$ then
 $C_s$ is closed, unbounded in $\mu^0_s,$ ordertype
$(C_s)\le\alpha(s).$ 

If $|s|$ is a successor ordinal then ordertype $(C_s)=\omega.$ 

(b) \ $\nu\in\Lim(C_s)\longrightarrow$ for some $\eta<|s|,$
$\nu=\mu^0_{s\restriction\eta }$ and $C_s\cap\nu=C_{s\restriction\eta }.$

(c) \ Let $\pi:\
\langle\overline{\Cal{A}},\overline{C}\rangle\overset{\Sigma_1}\to
\longrightarrow\langle {\Cal{A}}^0_s,C_s\rangle$ and write
crit$(\pi)=\bar\alpha,$ $\overline{\Cal{A}}=L_{\bar\mu }[\overline{A},\bar
s^*].$ If $\pi(\bar\alpha)=\alpha(s)$
  then $L[\overline{A},\bar s^*]\vDash|\bar s|$ is not a cardinal $>\bar
\alpha$  and 

(c1) \ $\overline{C}\in L_\mu [\overline{A},\bar s^*]$  where $\mu$ is the
least p\.r\. closed ordinal greater than $\bar\mu$ s\.t\. $L_\mu
[\overline{A},\bar s^*]\vDash$ card$(|\bar s|)\le\bar\alpha.$

(c2) \ $\pi$ extends to $\pi':\ {\Cal{A}}'\overset{\Sigma_1}\to\longrightarrow
{\Cal{A}}'_s$  where ${\Cal{A}}'=L_{\mu'}[\overline{A},\bar s^*],$ $\mu'=$
largest ordinal either equal to $\bar\mu$ or s\.t\. $\omega\cdot\mu'=\mu'$ and
$L_{\mu'}[\overline{A},\bar s^*]\vDash|\bar s|$ is a cardinal greater than
$\bar\alpha.$

(c3) \ If $\bar\alpha$ is a cardinal and $\pi(\bar\alpha)=\alpha$ then
$\overline{{\Cal{A}}}={\Cal{A}}^0_{\bar s}$  and $\overline{C}=C_{\bar s}.$

\vskip20pt

\flushpar
{\bf Relativized $\bold{\diamond}$ } \  Let $E=$ all $s\in S$ such that
$|s|$ limit and
ordertype $(C_s)=\omega.$  There exists $\langle D_s|s\in E\rangle$ such that
$D_s\subseteq {\Cal{A}}^0_s$ and:

(a) \  $D\in\widehat{{\Cal{A}}}_s\ne {\Cal{A}}^0_s,$ $D\subseteq
{\Cal{A}}^0_s\longrightarrow\{\xi<|s|\bigg|s\restriction\xi\in E,
D_{s\restriction\xi }=D\cap{\Cal{A}}^0_{s\restriction\xi }\}$ is stationary in
$\widehat{\Cal{A}}_s.$

(b) \ $D_s$ is uniformly $\Sigma_1$-definable as an element of ${\Cal{A}}'_s.$

(c) \ If ${\Cal{A}}'_s\vDash\alpha^{++}$ exists then $D_s=\emptyset.$

Now we use these combinatorial structures to impose some further restrictions
on membership in ${\Cal{P}}^u-{\Cal{P}}^{<u}.$   First some definitions. For
$p\in{\Cal{P}}^u$  and $\beta\in\Card\cap\alpha,$ $(p)_\beta$ denotes
$p\restriction\Card\cap[\beta,\alpha),$ $D\subseteq {\Cal{P}}^{<u}$ is {\it
predense} if every $p\in{\Cal{P}}^{<u}$ is compatible with an element of $D$
and for $\beta\in\Card\cap\alpha,D$ is {\it $\beta$-predense} if every
condition $p\in{\Cal{P}}^{<u}$ can be extended to some $q$  such that
$p\restriction\beta=q\restriction\beta$  and $q$ {\it meets} $D$  (i\.e\., $q$
extends an element of $D).$  And $p$  {\it reduces $D$ below $\beta$ } if
every $q\le p$ can be further extended to $r$  such that $r$  meets $D$  and
$(q)_\beta=(r)_\beta.$

\vskip20pt

\flushpar
{\bf Requirement A.} \ (Predensity Reduction) \  Suppose $p\in
{\Cal{P}}^u-{\Cal{P}}^{<u}.$

(A1) \  If $u\in E$ and $D_u\subseteq {\Cal{P}}^{<u}$ is $\beta$-predense
for all $\beta\in\Card\cap\alpha$ then $p$  meets $D_u.$

(A2) \ If $|u|$ is a successor ordinal, $D\subseteq {\Cal{P}}^{<u}$ is
predense and $D\in{\Cal{A}}^0_u$ then $p$  reduces $D$  below some
$\beta<\alpha.$

\vskip20pt

\flushpar
{\bf Requirement B.} \ (Restriction) \  For $p\in{\Cal{P}}^u$ let $|p|$ denote
the least $\xi$ s\.t\. $p\in {\Cal{P}}^{u\restriction\xi }.$  If $p$
belongs to ${\Cal{P}}^u$  and $\xi<|p|$  then there exists $r$ s\.t\.
$p\le r$  and $|r|=\xi.$ 

\vskip20pt

\flushpar
{\bf Requirement C.} \ (Nonstationary Restraint) \  Suppose
${\Cal{A}}_u\vDash\alpha$ inaccessible and $p\in{\Cal{P}}^u.$  Then there
exists a CUB $C\subseteq\alpha$ s\.t\. $C\in{\Cal{A}}_u$  and $\beta\in
C\longrightarrow\overline{p}_\beta=\emptyset.$ 

The remaining Requirement D will be introduced at a later point when we
discuss strong extendibility at successor stages.

Extendibility and distributivity for ${\Cal{P}}^u$ are stated as follows. Let
$q\le_\beta p$ signify that $q\le p$ and
$q\restriction\beta=p\restriction\beta.$  $({\Cal{P}}^{<u})_\beta$ denotes
$\{(p)_\beta|p\in {\Cal{P}}^{<u}\},$ for $\beta\in\Card\cap\alpha.$ 
$\underline{\Delta-}$ \underbar{distributivity} for ${\Cal{P}}^{<u}$ asserts
that if $D_\beta$  is $\beta^+$-predense on ${\Cal{P}}^{<u}$ for each
$\beta\in\Card\cap\alpha$ then every $p\in {\Cal{P}}^{<u}$ can be extended to
meet each $D_\beta.$ 

$(*)_u$ $p\in{\Cal{P}}^u,$ $\beta\in\Card\cap\alpha\longrightarrow\exists
q\le_\beta p$ $(q\in{\Cal{P}}^u-{\Cal{P}}^{<u})$

$(**)_u$ $({\Cal{P}}^{<u})_\beta$ is $\le\beta$-distributive in
$\widehat{\Cal{A}}_u$  for $\beta\in\Card\cap\alpha.$

And if $\alpha$ is inaccessible in ${\Cal{A}}^0_u$ then ${\Cal{P}}^{<u}$ is
$\Delta$-distributive in $\widehat{\Cal{A}}_u.$ 

\flushpar
These are proved by a simultaneous induction on $|u|.$  As the base case
$|u|=\alpha$  is vacuous we assume from now on that $|u|>\alpha.$  
The following
consequences of predensity reduction are needed in the proof.

\proclaim{Lemma 2.1} (Chain Condition for ${\Cal{P}}^{<u})$ \ Suppose $(**)_u$
holds. Then ${\Cal{P}}^{<u}$  has the $\alpha^+$-CC in $\widehat{\Cal{A}}_u.$ 
\endproclaim

\demo{Proof} We may assume that $\widehat{\Cal{A}}_u\ne {\Cal{A}}^0_u.$
Suppose $D\subseteq {\Cal{P}}^{<u}$ is predense and $D\in\widehat{\Cal{A}}_u.$
Consider $D^*=\{p\in{\Cal{P}}^{<u}|p$ reduces $D$ below some
$\beta\in\Card\cap\alpha\}.$  Then $D^*\in\widehat{\Cal{A}}_u.$  By $(**)_u$ 
and Lemma 1.2, $D^*$ is $\beta$-predense for all $\beta\in\Card\cap\alpha.$
(Use $\le\beta^+$-distributivity of $({\Cal{P}}^{<u})_{\beta^{+}}$ and
$\beta^{++}$-CC of $R^{G_{\beta^{+}}}\subseteq\beta^{++}$ denoting the
$({\Cal{P}}^{<u})_{\beta^{+}}$-generic, to reduce $D$  below $\beta^+.)$
Apply relativized $\diamond$ to obtain $\xi<|u|$ such that
$u\restriction\xi\in E,$
$D_{u\restriction\xi}=D^*\cap{\Cal{A}}^0_{u\restriction\xi}$  and
$D_{u\restriction\xi}$ is $\beta$-predense for all $\beta\in\Card\cap\alpha.$
Thus by predensity reduction and restriction,
$D^*\cap{\Cal{A}}^0_{u\restriction\xi }$  is predense on ${\Cal{P}}^{<u}$ and
therefore so is $D\cap{\Cal{A}}^0_{u\restriction\xi },$ a subset of $D$  of
$\widehat{\Cal{A}}_u$-cardinality $\le\alpha.$ \hfill{$\dashv$ }

\enddemo

\proclaim{Lemma 2.2} (Persistence for ${\Cal{P}}^{<u})$ \ Suppose $(**)_u$
holds, $D\subseteq {\Cal{P}}^{<u}$ is predense, $D\in\widehat{\Cal{A}}_u$  and
$u\subseteq v\in S_\alpha.$  Then $D$ is predense on ${\Cal{P}}^v.$ 
\endproclaim

\demo{Proof} By restriction, if $p\in{\Cal{P}}^v-{\Cal{P}}^u$ then $p$
extends some $q$ in ${\Cal{P}}^u-{\Cal{P}}^{<u}.$  By the chain condition for
${\Cal{P}}^{<u}$  we can assume that $D\in{\Cal{A}}^0_u$ and hence by
induction we can assume that $|u|$ is a successor ordinal. But then by
predensity reduction, $q$  reduces $D$  below some $\beta\in\Card\cap\alpha$
and hence so does $p.$  In particular $p$ is compatible with an element of
$D.$  \hfill{$\dashv$ }

\enddemo

We can now turn to the proofs of $(*)_u,$  $(**)_u.$  

\proclaim{Lemma 2.3} Assume $(**)_u$  and $|u|$  a limit ordinal. Then $(*)_u$
holds. 
\endproclaim

\demo{Proof} We first claim that if $p\in {\Cal{P}}^{<u}$  and $\langle
D_\beta|\beta_0\le\beta<\alpha\rangle$ $\in{\Cal{A}}^0_u,$
$D_\beta\subseteq{\Cal{P}}^{<u}$  $\beta^+$-predense for each $\beta$  then
there is $q\le_{\beta_{0}}p$ meeting each $D_\beta.$  We prove this with $
\alpha$  replaced by $\beta_1\in\Card\cap\alpha^+,$  by induction on
$\beta_1.$
The base case  $\beta_1=\beta^+_0$  and the case of $\beta_1$  a successor
cardinal  follow easily, using $(**)_u.$  If $\beta_1$  is singular in
${\Cal{A}}^0_u$ then we can choose $\gamma_0<\gamma_1<\cdots$ approximating $
\beta_1$  in length $\lambda<\beta_1$ and consider 
$\langle E_\delta|\delta<\lambda\rangle$ where $E_
\delta=$ all $q$ meeting each $D_\beta,$ $\lambda \le\beta<\gamma_\delta,|q|$
least so that $\langle
D_\beta|\beta_0\le\beta<\beta_1\rangle\in{\Cal{A}}^0_{u\restriction|q|}.$
Then
we are done by induction. If $\beta_1$ is inaccessible in ${\Cal{A}}^0_u$
then either $\beta_1=\alpha,$  in which case the result follows directly from
the second statement of $(**)_u,$ or $\beta_1<\alpha,$  in which case we can
factor ${\Cal{P}}^{<u}$ as $({\Cal{P}}^{<u})_{\beta_1^{+}}*
{\Cal{P}}^{G_{\beta_{1}^{+}}}$ (where $G_{\beta^{+}_{1}}$ denotes 
$\bigcup\{p_{\beta^{+}_{1}}|p\in G\},$ $G$ the generic for ${\Cal{P}}^{<u}).$
Then choose $(q)_{\beta^{+}_{1}}\le(p)_{\beta^{+}_{1}}$ that reduces each
$D_\beta,\beta_0\le\beta<\beta_1$ below $\beta^+_1,$ using $(**)_u$ and the
$\beta^+_1$-CC of ${\Cal{P}}^{G_{\beta^{+}_{1}}}.$ By induction on $\alpha,$
we can extend $q$  to meet all the $D_\beta$'s.

Now write $C_u=\{\mu^0_{u\restriction\xi_i}|i<\lambda\}$ and choose a
successor cardinal $\beta_0<\alpha$  to be at least as large as $\lambda$  and
the $\beta$  given in the statement of $(*)_u,$  if $\lambda<\alpha.$  Now
inductively define a subsequence $\langle\eta_j|j<\lambda_0\rangle$ of 
$\langle\xi_i|i<\lambda\rangle$ and conditions $\langle
p_j|j<\lambda_0\rangle$  as follows. First suppose $\lambda<\alpha.$  Let $p$
denote the condition given in the statement of $(*)_u.$  Set $p_0=p,\eta_0=$
least $\xi_i$ s\.t\. $p\in{\Cal{P}}^{<u\restriction\xi_i};$  $p_{j+1}=$ least
$q\le_\beta p_j$ s\.t\. for all $\gamma,\beta_0\le\gamma<\alpha,$ $q$ meets
all $\gamma^+$-predense $D\subseteq {\Cal{P}}^{<u\restriction\eta_j},$ $D\in
M^{\eta_{j}}_{\gamma^{+}}=\Sigma_1$ Skolem hull of $\gamma^+\cup\{p,\alpha\}$
in $\langle {\Cal{A}}^0_{u\restriction\eta_j},C_{u\restriction\eta_j}
\rangle,$ $\eta_{j+1}=$ least $\xi_i$ s\.t\.
$p_{j+1}\in{\Cal{P}}^{<u\restriction\xi_i};$ $p_\delta=$ g\.l\.b\. $\langle
p_j |j<\delta\rangle,$  $\eta_\delta=\bigcup\{\eta_j|j<\delta\}$  for limit
$\delta\le\lambda_0.$  The ordinal $\lambda_0$ is determined by the
condition that $\eta_{\lambda_{0}}$ is equal to $|u|.$  If $\lambda=\alpha$
then the definition is the same, except in defining $p_{j+1}$  require
$p_{j+1}\le_{\beta\cup\aleph_{i+1}}p_j$  where 
$\eta_j=\xi_i$ and only require $p_{j+1}$  to meet $\gamma^+$-predense
$D$  as above for $\gamma$ between $\beta\cup\aleph_i$  and $\alpha.$ 

We must verify that $p_\delta$ as defined above is indeed a condition for
limit $\delta.$ (There is no problem at successor stages, using Lemma 2.2 and
the first paragraph of the present proof.) First we show that for
$\gamma\in\Card 
\cap\alpha,$ $p_{\delta_{\gamma }}$ is reshaped. We need only consider
$\gamma\ge\beta$  and in case $\lambda=\alpha$  we need only consider
$\gamma\ge\beta\cup\aleph_i$  where $\eta_\delta=\xi_i.$  By construction if $
\gamma\in M^{\eta_{\delta }}_\gamma=\Sigma_1$ Skolem hull of $\gamma\cup\{p,
\alpha\}$ in $\langle {\Cal{A}}^0_{u\restriction\eta_\delta },C_{u\restriction
\eta_\delta}\rangle$  then $p_{\delta_{\gamma }}$ is
$\pi[({\Cal{P}}^{<u\restriction\eta_\delta })_\gamma ]$ generic over TC$(M^{
\eta_{\delta}}_\gamma)$  where $\pi:\ M^{\eta_{\delta }}_\gamma 
\longrightarrow TC(M^{\eta_{\delta }}_\gamma)$ is the transitive collapse. 
And $|p_{\delta_{\gamma }}|$ is $\Sigma_1$-definably singularized over TC$(M^{
\eta_{\delta }}_\gamma).$  Write TC$(M^{\eta_{\delta }}_\gamma)$ as 
$\langle\overline{\Cal{A}},\overline{C}\rangle.$ By genericity
and cofinality-preservation for
$\pi[({\Cal{P}}^{<u\restriction\eta}\delta)_\gamma ],$  
 $p_{
\delta_{\gamma }}$  codes $\overline{\Cal{A}}$  and by Relativized $\square$ 
(c1), $\overline{C}$ is constructible from $\overline{\Cal{A}}.$ So 
$p_{\delta_{\gamma }}$ is reshaped. If $M^{\eta_{\delta }}_{\gamma }
\cap\alpha =\gamma$ then  $p_{\delta_{\gamma }}$ is again reshaped
because of Relativized $\square$(c1), (no genericity argument
required). Lastly if $\gamma'=\min(M^{\eta_{\delta
}}_\gamma\cap(\ORD-\gamma))<\alpha$  then use the first argument, but with
$\pi[({\Cal{P}}^{<u\restriction\eta_\delta})_{\gamma'}]$ replacing 
$\pi[({\Cal{P}}^{<u\restriction\eta_{\delta }})_\gamma ].$ 

Next we show that $p_\delta\restriction\gamma\in
{\Cal{A}}_{p_{\delta_{\gamma }}}.$   As $p_\delta\restriction\gamma$ is
definable over TC$(M^{\eta_{\delta }}_\gamma)\in L[A\cap\gamma,
p_{\delta_{\gamma }}]$ this amounts to showing that $\mu_{p_{\delta_{\gamma
}}}$  is large enough. By $(**)_{u\restriction\eta_\delta }$ and Lemma 2.1 we
know that ${\Cal{P}}^{<u\restriction\eta_{\delta }}$ has the $\alpha^+$-CC
in $\widehat{\Cal{A}}_{u\restriction\eta_\delta }$  and hence
(when $M^{\eta_{\delta }}_\gamma\cap\alpha\ne\gamma)$  $p_{\delta_{\gamma }}$
is in fact $\pi^{'-1}[({\Cal{P}}^{<u\restriction\eta_\delta
})_{\gamma'}]$-generic over ${\Cal{A}}',$  where $\pi'$
(with domain ${\Cal{A}}')$   is the extension of $\pi^{-1}$ given by
Relativized $\square$ (c2) and $\gamma'=\min(M^{\eta _{\delta }}_\gamma
\cap(\ORD-\gamma )).$  And thus ${\Cal{A}}'[p_{\delta_{\gamma }}]\vDash|
p_{\delta_{\gamma }}|$ is a cardinal. But by Relativized $\square$ (c1), 
TC$(M^{\eta _{\delta }}_\gamma )$ appears relative to $p_{\delta_{\gamma }}$
before the next p\.r\. closed ordinal after the height of ${\Cal{A}}'.$  So
$p_\delta\restriction\gamma\in {\Cal{A}}_{p_{\delta_{\gamma }}}.$  If
$M^{\eta_{\delta }}_\gamma\cap\alpha=\gamma$  then no genericity argument is
required; we only need Relativized $\square$ (c1).

Requirements B, C  are easily checked, the latter using the fact that in case
of $\alpha$  inaccessible in ${\Cal{A}}^0_u$ we required
$p_{j+1}\le_{\beta\cup\aleph_{i+1}}p_j(\eta_j=\xi_i)$  and
therefore can use diagonal intersection of clubs. To check Requirement (A1)
note that if $M^{\eta_{\delta }}_\gamma\cap\alpha\ne\gamma$  then either
$p_{\delta_{\gamma }}\notin E$ or $D_{p_{\delta_{\gamma }}}=\emptyset,$ since
${\Cal{A}}'_{p_{\delta_{\gamma }}}\vDash\gamma^{++}$ exists and we can apply 
Relativized $\diamond$ (c). If $M^{\eta_{\delta }}_\gamma\cap\alpha=\gamma$
then $p_{\delta_{\gamma }}\in E$ iff $u\restriction\eta_\delta\in E$  by
Relativized $\square$  (c3) and if these hold then by Relativized $\diamond$
(b), $\pi'[D_{p_{\delta_{\gamma }}}]=D_{u\restriction\eta_\delta },$ where
$\pi'$  comes from Relativized $\square$ (c2). So all we need to arrange is
that our initial condition $p$  be chosen to meet $D_u,$ in case $u\in E,$
and otherwise choose $\eta_0$ to be at least $\xi_\omega,$  so that
$u\restriction\eta_\delta\notin E$  for limit $\delta.$ \hfill{$\dashv$}

\enddemo

\proclaim{Lemma 2.4} Assume $|u|$ limit and $(*)_v,(**)_v$ for $v\subseteq
u,v\ne u.$  Then $(**)_u$ holds.
\endproclaim

\demo{Proof} We may assume that $\widehat{\Cal{A}}_u\ne{\Cal{A}}^0_u.$  We
need only make a small change in the construction of the proof of Lemma 2.3.
Given predense $\langle D_i|i<\beta\rangle$ on $({\Cal{P}}^{<u})_\beta$ 
in $\widehat{\Cal{A}}_u$ 
with 
$\beta<\alpha,$  select $\xi<|u|$ of cofinality $>\beta$ such that $D_i\cap
({\Cal{P}}^{<u\restriction\xi })_\beta$ is predense on
$({\Cal{P}}^{<u\restriction\xi })_\beta$  for all $i<\beta$  and then choose
the continuous sequence $\langle\xi_i|i<\beta\rangle$  from
$C_{u\restriction\xi }$ by: $\xi_0=\omega ^{\text{th}}$ element of 
$C_{u\restriction\xi },$ $\xi_{i+1}=$
least $\xi^*\in C_{u\restriction\xi }$ greater than $\xi_i$ s\.t\.
$q\in({\Cal{P}}^{u\restriction\xi_i})_\beta\longrightarrow\exists r\le
q(r\in({\Cal{P}}^{u\restriction\xi^*})_\beta,r$ meets $D_i),$
$\xi_\lambda=\bigcup\{\xi_i|i<\lambda\}$ for limit $\lambda\le\beta.$  Then
$u\restriction\xi_\lambda\notin E$ and $\langle\xi_i|i<\lambda\rangle\in
{\Cal{A}}_{u\restriction\xi_\lambda }$ for limit $\lambda.$ 
\enddemo

Now repeat the construction of the proof of Lemma 2.3, extending along the
$\xi_i$'s instead of along $C_u,$ hitting $D_i$  at stage $i+1.$ We can
guarantee  $\langle
D_i\cap({\Cal{P}}^{<u\restriction\xi_i})_\beta|i<\lambda\rangle$ is
$\Delta_1\langle {\Cal{A}}^0_{u\restriction\xi_\lambda
},C_{u\restriction\xi_\lambda}\rangle$  in our choice of $\xi_i$'s
as well, so hitting the $D_i$'s does not interfere with the proof that
$p_\delta$  is a condition for limit $\delta.$  The proof of
$\Delta$-distributivity is similar. \hfill{$\dashv$ }

\proclaim{Lemma 2.5} Suppose $(**)_u$ holds and $|u|$ is a successor ordinal.
Then $(*)_u$ holds.
\endproclaim

\demo{Proof} We may assume that the given $p$  belongs to
${\Cal{A}}_v-{\Cal{A}}^0_v$ where $v=u\restriction(|u|-1).$  Write $C_u=
\langle\xi_j|j<\omega\rangle.$ Now proceed as in the construction of the proof
of Lemma 2.3, making successive $\le_\beta$-extensions below $p$ (where
$\beta$ is given in the statement of $(*)_u),$  $p\ge_\beta p_0\ge_\beta
p_1\ge_\beta\cdots$  so that $p_{j+1}$ meets all $\gamma^+$-predense 
$D\subseteq{\Cal{P}}^{<u}$ in $M^{\xi_{j}}_{\gamma^{+}},$ where
$M^{\xi_{j}}_{\gamma^{+}}=\Sigma_1$ Skolem hull of
$\gamma^+\cup\{p,\alpha,\xi_0,\cdots,\xi_{j-1}\}$  in
${\Cal{A}}_v\restriction\xi_j,$  for all $\gamma\in [\beta,\alpha).$  If we
set $\hat q=$ g\.l\.b\. $\langle p_i|i\in\omega\rangle$ then $\hat q$ meets
the requirements for being a condition at all $\gamma\in\Card\cap\alpha^+$
with the exception of $\gamma$ in $C\cup\{\alpha\},$
$C=\{\gamma|M_\gamma \cap\alpha=\gamma\},$  $M_\gamma=\Sigma_1$ Skolem hull of
$\gamma\cup\{p,\alpha\}$ in $\langle {\Cal{A}}_v,C_u\rangle.$  The reason is
that for $\gamma\in\alpha-C,T_\gamma=TC(M_\gamma)$ belongs to ${\Cal{A}}_{\hat
q_\gamma },$  since $T_\gamma\vDash|\hat q_\gamma|$ is a cardinal and $\hat
q_\gamma$ is generic over $T_\gamma.$ 

To make $\hat q$ into a condition $q\in {\Cal{P}}^u$ we must do two things.
First extend $\hat q_{\gamma^{+}}$ for $\gamma\ge\beta$ so as to code
$u(|v|)=0$ or $1.$  This is easily done as there are no conflicts between the
successor and limit codings. Second for $\gamma\in C$  we extend $\hat
q_\gamma$  to $q_\gamma=\hat q_\gamma *u(|v|).$  The only remaining question
is whether the reatraint $\overline{\hat q}_\gamma$ will allow us to do this.
But $\gamma\in C\longrightarrow\overline{\hat q}_\gamma=\emptyset$ since $C$
is contained in the CUB witnessing Requirement C for $\hat q$ at $\alpha.$ 
\hfill{$\dashv$ }
\enddemo

\proclaim{Lemma 2.6} Suppose $(*)_u$ and $(**)_v,v\subseteq u\ne v$ hold and
$|u|$ is a successor. Then $(**)_u$ holds.
\endproclaim

\demo{Proof} We must show that if $v=u\restriction(|u|-1)$  and
$p\in({\Cal{P}}^v)_\beta-({\Cal{P}}^{<v})_\beta,$ $\langle
D_i|i<\beta\rangle\in {\Cal{A}}_v$ are predense on $({\Cal{P}}^v)_\beta$  then
there exists $q\le p$ meeting each $D_i.$  For simplicity we assume
$\beta=\omega.$  

\flushpar
{\bf Definition.} \ Suppose  $f(\beta)=M_\beta$ is a function in ${\Cal{A}}_v$
from $\Card^+\cap\alpha$ $(\Card^+$ denotes all successor cardinals) into
${\Cal{A}}_v$ such that card $(M_\beta)\le\beta$  for all $\beta\in\Dom(f)$
and suppose $p\in {\Cal{P}}^v.$  Then $\Sigma_f^p=\{q\in{\Cal{P}}^v|\ \forall\
\beta\in\Dom(f),$ $q(\beta)$  meets all predense $D\subseteq
R^{p_{\beta^{+}}},$  $D\in M_\beta\}.$ 
\enddemo

\proclaim{Sublemma 2.7} $\Sigma^p_f$ is dense below $p$  in ${\Cal{P}}^v.$ 
\endproclaim

Before proving Sublemma 2.7 we establish the Lemma, assuming it. Choose a
limit ordinal $\lambda=\omega^\lambda<\mu_v$ such that $\langle
D_i|i<\omega\rangle,$ $C_v\in{\Cal{A}}_v\restriction\lambda$  $=L_\lambda [A\cap\alpha,v^*]$ 
and
$\Sigma_1\cof({\Cal{A}}_v\restriction\lambda)=\omega.$  Choose a
$\Sigma_1({\Cal{A}}_v\restriction\lambda)$ sequence
$\lambda_0<\lambda_1<\cdots$  cofinal in $\lambda$ such that $\langle
D_i|i<\omega\rangle,C_v,$ $x\in {\Cal{A}}_v\restriction\lambda_0$ where $x$ is
a parameter defining the $\lambda_i$'s. Set $M^i_\gamma=$ least
$M\prec_{\Sigma_1}A_v\restriction\lambda_i$ such that $\gamma\cup\{x,\langle
D_i|i<\omega\rangle,\alpha,C_v\}\subseteq M,$ for each
$\gamma\in\Card^+\cap\alpha.$ Define $f_i(\gamma)=M^i_\gamma.$ 

Choose $p=p_0\ge p_1\ge\cdots$ successively so that $p_{i+1}$ meets $D_i$ and
$\Sigma_{f_{i}}^{p_{i}}.$ Set $p^*=$ g\.l\.b\. $\langle
p_i|i\in\omega\rangle.$ We show that $p^*$ is a well-defined condition.
If $|v|>\alpha$ then 
thanks to $(**)_v$ it will suffice to show that if $D\in
M^i_\gamma\cap{\Cal{A}}^0_v$  is predense on $({\Cal{P}}^{<v})_\gamma,$  
$\gamma\in\Card\cap\alpha$ then some $p_j$  reduces $D$  below $\gamma.$  (For
then, $p^*_\gamma$ codes a generic over the transitive collapse of
$M^i_\gamma\cap{\Cal{A}}^0_v.)$  If $|v|=\alpha$ then instead of
${\Cal{P}}^{<v}=\emptyset$ use ${\Cal{P}}^\alpha=\{p\restriction\beta^+\big|
\beta\in\Card\cap\alpha,\ p\in {\Cal{P}}^v\},$ ordered in the natural way.
Note that ${\Cal{P}}^\alpha$ is cofinality-preserving, by applying $(**)$ at
cardinals $<\alpha.$

Choose $j\ge i$ so that for $k>j,p_k$  reduces $D$  no further than $p_j.$
Let $\gamma'$  be least so that $p_j$  reduces $D$  below $\gamma'.$  Then
$\gamma'<\alpha$  by Predensity Reduction for $p.$  If $\gamma'\le\gamma$ then
of course we are done. If $\gamma'>\gamma$ is a double successor cardinal then
we reach a contradiction since by definition $p_{j+1}$ reduces $D$  further.
If $\gamma'=\delta^+,$ $\delta$ a limit cardinal then by
Predensity Reduction at $\delta,D$ is reduced below some
$\delta'<\delta,$ another contradiction. If $\gamma'$  is a limit cardinal
then the same argument applies, replacing $\gamma'$ by $(\gamma')^+.$  

Finally we have:

\demo{Proof of Sublemma 2.7} It suffices to show the following.

\flushpar
{\bf Strong Extendibility}  \ Suppose $g\in {\Cal{A}}_v,g(\beta)\in
H_{\beta^{++}}$  for all $\beta\in\Card\cap(\beta_0,\alpha)$  and $p\in
{\Cal{P}}^v.$  Then there is $q\le_{\beta_{0}}p$ such that
$g\restriction\beta\in{\Cal{A}}_{q_{\beta }}$  for all
$\beta\in\Card\cap(\beta_0,\alpha ].$  

For, Strong Extendibility allows us to extend to a condition $q$  such that
for all $\beta\in\Card\cap\alpha,$ $g\restriction\beta\in{\Cal{A}}_{q_{\beta
}},$  where $g(\beta)=f(\beta)\cap H_{\beta^{++}}.$  Then successively
extend each $q(\beta)$ to meet predense $D$ in $f(\beta).$  

We now break down Strong Extendibility into the ramified form in which it will
be proved. For any $\mu$  such that $\mu^0_v\le\mu<\mu_v,$ $k\in\omega-\{0\}$
and $\beta\in\Card\cap\alpha$ let $M^{\mu,k}_\beta=\Sigma _k$ Skolem  hull
of $\beta\cup\{\alpha\}$ in ${\Cal{A}}^*_v\restriction\mu=\langle L_\mu
[A\cap\alpha,v^*,C_v],$  $A\cap\alpha,v^*,C_v\rangle.$  (Notice that this
structure is $\Sigma_1$ projectible to $\alpha$  without parameter.)

\flushpar
${\underline{SE(\mu ,k)}}$ \ Suppose $p\in{\Cal{P}}^v$ and
$\beta_0\in\Card\cap\alpha.$  Then there exists $q\le_{\beta_{0}}p$  such that
$TC(M^{\mu,k}_\beta)\in{\Cal{A}}_{q_{\beta }}$  for all
$\beta\in\Card\cap(\beta_0,\alpha).$  

It suffices to prove $SE(\mu,k)$ for all  $\mu,k$ as above.  We do so by
induction on $\mu$  and for fixed $\mu,$  by induction on $k.$  To verify
the base case of this
induction  we must impose one last requirement on our
conditions. 

\flushpar
{\bf Requirement D} \ Suppose $p\in {\Cal{P}}^v-{\Cal{P}}^{<v}$  and  $g\in
{\Cal{A}}^0_v,$ $g(\beta)\in H_{\beta^{++}}$  for all
$\beta\in\Card\cap\alpha.$  Then $g\restriction\beta\in{\Cal{A}}_{p_{\beta }}$
for sufficiently
large $\beta\in\Card\cap\alpha.$ 

This requirement is respected by our earlier constructions. Now, if
$k=1$  and
$\mu$  is a limit ordinal then we can use a
$\Sigma_1({\Cal{A}}^*_v\restriction\mu)$  approximation to $\mu$  and
induction (or Requirement $D$ if $\mu =\mu ^0_v)$ 
to obtain 
$q\le p$  satisfying the conclusion of $SE(\mu,1),$ using
the $\Sigma_f$'s for $f\in {\Cal{A}}^*_v\restriction\mu.$  Similarly if $\mu$
is a successor, $k=1$ then use
$\langle\Sigma_k({\Cal{A}}^*_v\restriction\mu-1)|k\in\omega\rangle$  to
approximate $\Sigma_1({\Cal{A}}^*_v\restriction\mu),$ using the $\Sigma_f$'s,
$f$  definable over ${\Cal{A}}^*_v\restriction\mu-1.$ 

Suppose $k>1.$  By induction we can assume that TC$(M^{\mu,k-1}_\beta)\in
{\Cal{A}}_{p_{\beta }}$ for large enough $\beta.$  If
$C=\{\beta<\alpha|\beta=\alpha\cap M^{\mu,k}_\beta\}$  is unbounded in
$\alpha$  then successively extend $p\restriction\beta$ for $\beta\in C$ 
 so that
TC$(M^{\mu,k}_{\beta'})\in{\Cal{A}}_{q_{\beta'}}$  for $\beta'<\beta.$  There
is no problem at limits since TC$(M^{\mu,k}_\beta),$
$C\cap\beta\in{\Cal{A}}_{p_{\beta }}$  for $\beta\in C.$ 

If $\alpha$ is $\Sigma_k({\Cal{A}}^*_v\restriction\mu)$-singular then choose a
continuous cofinal $\Sigma_k({\Cal{A}}^*_v\restriction\mu)$  sequence
$\beta_0<\beta_1<\cdots$  below $\alpha$  of ordertype
$\lambda_0=\cof(\alpha).$  Also choose $\beta_{i+1}$ large enough so that
$M^{\mu,k-1}_{\beta_{i+1}}\models\beta_i$ is defined. This is possible since
${\Cal{A}}^*_v\restriction\mu=\bigcup\{M^{\mu,k-1}_\beta|\beta<\alpha\}.$  Now
define $N^i_\beta$ for $i<\lambda_0,$ $\beta<\beta_i$ to be the $\Sigma_k$ 
Skolem hull of $\beta$  in $M_{\beta_{i}}^{\mu,k-1}.$ Then $\langle
TC(N^i_\beta)|\beta<\beta_i\rangle\in{\Cal{A}}_{p_{\beta_{i}}}$ for
$i<\lambda_0$ since it is easily defined from
$M_{\beta_{i}}^{\mu,k-1}\in{\Cal{A}}_{p_{\beta_{i}}}.$  Successively
$\lambda_0$-extend $p\restriction\beta_i,$ producing
$p=p_0\ge_{\lambda_{0}}p_1\ge_{\lambda_{0}}\cdots$ where
$TC(N^i_\beta)\in{\Cal{A}}_{p_{i_{\beta }}}$  for
$\beta\in(\lambda_0,\beta_i).$  This is possible by induction on $\alpha,$
and since TC$(N^i_\beta)$ is easily defined from $\langle TC(N^i_{\bar\beta
})|\bar\beta<\beta\rangle$  for limit $\beta<\beta_i.$ (We must also require
that $p_{i+1}$ meets $\Sigma^{p_{i}}_{f_{i}}$ where $f_i(\beta)=N^i_\beta.)$ 
$p_\lambda$  is well-defined for limit $\lambda\le\lambda_0$  and
${\Cal{A}}_{p_{\lambda_{0_{\beta }}}}$ contains $\langle
TC(N^i_\beta)|i<\lambda_0\rangle$ and hence TC$(M_\beta^{\mu,k})$  for
$\beta>\lambda_0.$  Then use induction to fill in on $(0,\lambda_0]$ so
that $SE(\mu,k)$ is satisfied.

Lastly, there is the intermediate case where $\alpha$ is
$\Sigma_k({\Cal{A}}^*_v\restriction\mu)$-regular but 
$C=\{\beta <\alpha |\beta=\alpha\cap M_\beta^{\mu,k}\}$ is bounded in
$\alpha.$  Then 
$\Sigma_{k+1}({\Cal{A}}^*_v\restriction\mu)$-cof $(\alpha)=\omega$  and we
apply induction to produce $p=p_0\ge p_1\ge\cdots$ so that
$p_{i+1}\restriction [\beta_i,\beta_{i+1}]$  obeys $SE(\mu,k)$ where
$\beta_0<\beta_1<\cdots$ is a cofinal $\omega$-sequence of successor cardinals
below $\alpha.$  Let $q=$ g\.l\.b\. $\langle p_i|i\in\omega\rangle.$ 

This completes the proof of Sublemma 2.7 and hence of $(**)_u.$ 
\hfill{$\dashv$ }
\enddemo

\vskip20pt

\flushpar
{\bf Section Three \ Proof of Jensen's Theorem}

A condition in ${\Cal{P}}$ is a function $p$  from an initial segment of Card
into $V$  such that $\Dom(p)$ has a maximum $\alpha(p),$ for any
$\alpha\in\Dom(p),p(\alpha)=(p_\alpha,\bar p_\alpha),$ if
$\alpha\in\Dom(p)\cap\alpha(p)$  then $p(\alpha)$ belongs to
$R^{p_{\alpha^{+}}},$  $p(\alpha(p))=(s(p),\emptyset)$  where $s(p)\in
S_{\alpha(p)}$  and for uncountable limit cardinals $\alpha\in\Dom(p),$
$p\restriction\alpha\in{\Cal{P}}^{p_{\alpha }}.$  And $q\le p$ in ${\Cal{P}}$
if $\alpha(p)\le\alpha(q),$ $s(p)\subseteq q_{\alpha(p)}$ and for
$\alpha\in\Dom(p)\cap\alpha(p),$ $q(\alpha)\le p(\alpha)$  in
$R^{q_{\alpha^{+}}}.$  

For any $\alpha\in\Card, s\in S_\alpha,$ ${\Cal{P}}^s$ denotes all
$p\restriction\alpha$ for $p\in{\Cal{P}}$ such that $\alpha(p)=\alpha$  and
$s(p)=s.$  And ${\Cal{P}}^\alpha$  denotes all $p\in{\Cal{P}}$  such that
$\alpha(p)<\alpha.$ 

Now suppose $\alpha$  is an uncountable limit cardinal and $s\in S_\alpha,$
$|s|=\alpha+1.$  By Lemma 2.2, $G$ ${\Cal{P}}$-generic $\longrightarrow
G\cap{\Cal{P}}^{<s}$  is ${\Cal{P}}^{<s}$-generic over ${\Cal{A}}^0_s=L_\mu
[A\cap\alpha ],\mu$  the least p\.r\. closed ordinal greater than $\alpha.$
As the forcing relation for ${\Cal{P}}^{<s}$  restricted to sentences of rank
$<\alpha$  belongs to $L_\mu [A\cap\alpha ],$  it follows that the forcing
relation $p\Vdash\varphi,$ $p\in{\Cal{P}}$ and $\varphi$ ranked, is 
$\langle L[A],A\rangle$-definable: $p\Vdash\phi$ iff for some $\alpha$  as
above, $\phi$ has rank $<\alpha,p\in L_\alpha [A]$  and $L_\mu [A\cap\alpha
]\models$ ``$p\Vdash\phi$'', $\mu$ the least p\.r\. closed ordinal $>\alpha.$

Now note that ${\Cal{P}}$  preserves cofinalities, as otherwise ${\Cal{P}}^s$
would change cofinalities for some $s$  as above, contradicting Distributivity
(Lemmas 2.4, 2.6) and Chain Condition (Lemmas 1.2, 2.1). If $G$  is
${\Cal{P}}$-generic then $L[G]=L[X]$  where $X=G_\omega\subseteq\omega_1.$
Finally by Jensen-Solovay [68], $X$  can be coded by a real via a CCC forcing.
This completes the proof of Jensen's Coding Theorem, subject to the
verification of Relativized $\square$  and $\diamond.$ 

\vskip20pt

\flushpar
{\bf Section Four \  Relatived Square and Diamond}

For completeness, we prove Relativized $\square$ and $\diamond.$  As
relativization causes no serious problems, we first establish unrelativized
versions, and then afterward indicate what modifications are required. We
begin with $\square.$

First we prove $\square$ in the following form:

\flushpar
{\bf Global $\bold{\square}$ } \  Assume $V=L.$  Then there exists $\langle
C_\mu|\mu$  a singular limit ordinal$\rangle$ such that:

(a) \  $C_\mu$ is CUB in $\mu$

(b) ordertype $(C_\mu)<\mu$ 

(c) $\bar\mu\in\Lim C_\mu\longrightarrow C_{\bar\mu }=C_\cap\bar\mu.$ 

\flushpar
In the proof we shall take advantage of Jensen's $\Sigma^*$ theory, as
reformulated in Friedman [94]. For the convenience of the reader we describe
that theory here. 

For simplicity of notation, for limit ordinals $\mu$  we let 
$\widetilde J_\mu$
denote $J_\alpha$  where $\omega\alpha=\mu.$  So ORD$(\widetilde
J_\mu)=\mu.$

Let $M$  denote some $J_{\alpha },\alpha>0.$  (More generally, our theory
applies to ``acceptable $J$-models''.) We make the following definitions,
inductively. We order finite sets of ordinals by the maximum difference order:
$x<y$ iff $\alpha\in Y$  where $\alpha$  is the largest element of
$(y-x)\cup(x-y).$  

1) \ A $\Sigma^*_1$ formula is just a $\Sigma_1$ formula. A predicate is 
$\underline{\Sigma^*_1}$ ($\Sigma^*_1,$ respectively) if it is definable by a
$\Sigma^*_1$ formula with (without, respectively) parameters.
$\rho^M_1=\Sigma^*_1$  projectum of $M=$ least $\rho$  s\.t\. there is a
$\underline{\Sigma^*_1}$  subset of $\omega\rho$  not in $M$
and $p^M_1=$ least
$p$  s\.t\. $A\cap\rho_1^M\notin$ $M$ for some $A$ $\Sigma^*_1$ in parameter
$p$  (where $p$  is a finite set of ordinals).
$H^M_1=H^M_{\omega\rho_{1}^{M}}=$ sets $x$ in $M$  s\.t\. $M$-card (transitive
closure $(x))<\omega\rho^M_1.$  For any $x\in M,$ $M_1(x)=$ First reduct of
$M$  relative to $x=\langle H^M_1, A_1(x)\rangle$  where $A_1(x)\subseteq
H^M_1$  codes the $\Sigma^*_1$  theory of $M$  with parameters from
$H^M_1\cup\{x\}$  in the natural way: $A_1(x)=\{\langle y,n\rangle|$  the
$n^{\text{th}}$ $\Sigma^*_1$ formula is true at $\langle y,x\rangle,$ $y\in
H^M_1\}.$  A good $\Sigma^*_1$ function is just a $\Sigma_1$ function and for
any $X\subseteq M$  the $\Sigma_1^*$ hull $(X)$  is just the $\Sigma_1$ hull
of $X.$

(2) For $n\ge 1,$  a $\Sigma^*_{n+1}$  formula is one of the form
$\phi(x)\longleftrightarrow M_n(x)\models\psi,$ where $\psi$ is $\Sigma_1.$  A
predicate is $\underline{\Sigma^*_{n+1}}$ ($\Sigma^*_{n+1},$  respectively) if
it is defined by a $\Sigma^*_{n+1}$ formula with (without, respectively)
parameters. $\rho^M_{n+1}=\Sigma^*_{n+1}$ projectum of $M=$ least $\rho$  such
that there is a $\underline{\Sigma^*_{n+1}}$ subset of $\omega\rho$ not in
$M$ 
and $p^M_{n+1}=p^M_n\cup p$ where
$p$ is least such that $A\cap\rho_{n+1}^M\notin$ $M$ for some $A$ 
$\Sigma^*_{n+1}$ in parameter
$p^M_n\cup p.$
  $H^M_{n+1}=H^M_{\omega\rho_{n+1}^{M}}=$ sets $x$ in $M$  s\.t\.
$M$-card (transitive closure $(x)$) $<\omega\rho^M_{n+1}.$  For any $x\in M,$
$M_{n+1}(x)=(n+1)$ st reduct of $M$  relative to $x=\langle H^M_{n+1},
A_{n+1}(x)\rangle$  where $A_{n+1}(x)\subseteq H^M_{n+1}$ codes the
$\Sigma^*_{n+1}$  theory of $M$  with parameters from $H^M_{n+1}\cup\{x\}$ in
the natural way: \ $A_{n+1}(x)=\{\langle y,m\rangle|$  the $m^{\text{th}}$ 
$\Sigma^*_{n+1}$ formula is true at $\langle y,x\rangle,y\in H^M_{n+1}\}.$  A
good $\Sigma^*_{n+1}$  function $f$  is a function whose graph is
$\Sigma^*_{n+1}$  with the additional property that for $x\in\Dom(f),$
$f(x)\in\Sigma^*_n$ hull $(H^M_n\cup\{x\}).$  The $\Sigma^*_{n+1}$ hull $(X)$
for $X\subseteq M$ is the closure of $X$  under good $\Sigma^*_{n+1}$
functions. 

\vskip20pt

\flushpar
{\bf Facts.} \ (a) \   $\phi,\psi\Sigma^*_n$  formulas
$\longrightarrow\phi\vee\psi,\phi\wedge\psi$  are $\Sigma^*_n$  formulas

(b) \   $\phi\Sigma^*_n$ or $\prod^*_n$ ($=$ negation of
$\Sigma^*_n)\longrightarrow\phi$ is $\Sigma^*_{n+1}$

(c)  \  $Y\subseteq\Sigma^*_n$ hull $(X)\longrightarrow\Sigma^*_n$ hull
$(Y)\subseteq\Sigma^*_n$ hull $(X)$

(d) \ $f$  good $\Sigma^*_n$ function $\longrightarrow f$ good
$\Sigma^*_{n+1}$ function

(e) \ $\Sigma^*_n$ hull $(X)\subseteq\Sigma^*_{n+1}$ hull $(X)$

(f) \  There is a $\Sigma^*_n$ relation $W(e,x)$ s\.t\. if $S(x)$  is
$\Sigma^*_n$ then for some $e\in\omega,$ $S(x)\longleftrightarrow W(e,x)$  for
all $x.$

(g) \ The structure $M_n(x)=\langle H^M_n, A_n(x)\rangle$ is amenable.

(h) \ $H^M_n=J^{A_{n}}_{\omega\rho^{M}_{n}}$  where $A_n=A_n(0).$

(i) \ Suppose $H\subseteq M$ is closed under good $\Sigma^*_n$ functions and
$\pi:\ \overline{M}\longrightarrow M,$ $\overline{M}$ transitive, Range
$(\pi)=H$  and $p^M_{n-1}\in H$ (if $n>1).$ 
  Then $\pi$ preserves $\Sigma^*_n$ formulas: for $\Sigma^*_n\phi $
and $x\in\overline{M},$ $\overline{M}\models\phi(x)\longleftrightarrow
M\models\phi(\pi(x)).$ And (for $n>1),$ $\pi(p^{\bar M}_{n-1})=p_{n-1}^M.$    

\flushpar
{\bf Proof of (i)} \ Note that   $H\cap M_{n-1}(\pi(x))$ is $\Sigma
_1$-elementary in $M_{n-1}(\pi(x)).$  And $\pi^{-1}[H\cap
M_{n-1}(\pi(x))]=\langle J^A_{\omega\rho },A(x)\rangle$  for some
$\rho,A,A(x).$  But (by induction on $n)$  $A=A^M_{n-1}\cap J^A_{\omega\rho
},$ 
$A(x)=A_{n-1}(x)^{\bar M}\cap J^A_{\omega\rho }.$  And $\rho=\rho^{\bar M}_{n-1}$
using our assumption about the parameter $p^M_{n-1}.$  And
$\pi^{-1}(p^M_{n-1})=\bar p$  must be $p^{\bar M}_{n-1}$ as $\bar
M=\Sigma^*_{n-1}$ hull of $H^{\bar M}_{n-1}\cup\{p_{n-1}^{\bar M}\}.$  
\hfill{$\dashv$ }

\proclaim{Theorem 4.1} By induction on $n>0:$

1) \ If $\phi(x,y)$ is $\Sigma^*_n$ then $\exists
y\in\Sigma^*_{n-1}$  hull $(H^M_{n-1}\cup\{x\})\phi(x,y)$ is also
$\Sigma^*_n.$

2) \ If $\phi(x_1\cdots x_k)$ is $\Sigma^*_m,m\ge n$  and $f_1(x),\cdots,
f_k(x)$  are good $\Sigma^*_n$  functions, then $\phi(f_1(x)\cdots f_k(x))$ is
$\Sigma^*_m.$

3) \ The domain of a good $\Sigma^*_n$ function is $\Sigma^*_n$

4) \ Good $\Sigma^*_n$ functions are closed under composition.

5) \ ($\Sigma^*_n$ Uniformization) \ If $R(x,y)$ is $\Sigma^*_n$  then there
is a good $\Sigma^*_n$  function $f(x)$ s\.t\. $x\in\Dom(f)\longleftrightarrow
\exists y\in\Sigma^*_{n-1}$ hull
$(H^M_{n-1}\cup\{x\})R(x,y)\longleftrightarrow R(x,f(x)).$  

6) \ There is a good $\Sigma^*_n$ function $h_n(e,x)$ s\.t\. for each
$x,\Sigma^*_n$ hull $(\{x\})=\{h_n(e,x)|e\in\omega\}.$ 
\endproclaim

\demo{Proof}  The base case $n=1$ is easy (take $\so$ hull $(X)=M$ for all
$X).$   Now we prove it for $n>1,$ assuming the result for smaller $n.$

1) \  Write $\exists y\in\sgn$ hull $(H^M_{n-1}\cup\{x\})\phi(x,y)$ as 
$\exists\bar y\in H^M_{n-1}\phi(x,h_{n-1}(e,\langle x,\bar y\rangle))$ using
6) for $n-1.$  Since $h_{n-1}$ is good $\sgn$ we can apply 2) for $n-1$ to
conclude that $\phi(x,h_{n-1}(e,\langle x,\bar y\rangle))$ is $\sg.$ Since the
quantifiers $\exists e\exists\bar y\in H^M_{n-1}$ range over $H^M_{n-1}$ they
preserve $\sg$-ness.

2)  \  $\phi(f_1(x)\cdots f_k(x))\longleftrightarrow\exists x_1\cdots
x_k\in\sgn$ hull $(H^M_{n-1}\cup\{x\})$  $[x_i=f_i(x)$ for $1\le i\le
k\wedge\phi(x_1\cdots x_k)].$  If $m=n$ then this is $\sg$ by 1). If $m>n$
then reason as follows: the result for $m=n$  implies that $A_n(\langle
f_1(x)\cdots f_k(x)\rangle)$  is $\Delta_1$ over $M_{n+1}(x).$  Thus
$A_{m-1}(\langle f_1(x)\cdots f_k(x)\rangle)$ is $\Delta_1$ over $M_{m-1}(x).$
So as $\phi$ is $\Sigma^*_m$ we get that $\phi(f_1(x)\cdots f_k(x))$ is also
$\Sigma_1$  over $M_{m-1}(x),$  hence $\Sigma^*_m.$

3) \  If $f(x)$ is good $\Sigma^*_n$ then dom$(f)=\{x|\exists
y\in\Sigma^*_{n-1}$ hull of $H^M_{n-1}\cup\{x\}(y=f(x))\}$ is $\sg$ by  1).

4)  \  If $f,g$ are good $\sg$ then the graph of $f\circ g$ is $\sg$  by 2).
And $f\circ g(x)\in\sgn$ hull$(H^M_{n-1}\cup\{x\})$ since the latter hull
contains $g(x),f$ is good $\sg$ and Fact c) holds.

5)  \  Using 6) for $n-1,$ let $\overline{R}(x,\bar y)\longleftrightarrow
R(x,h_{n-1}(\bar y))\wedge\bar y\in H^M_{n-1}.$  Then $\overline{R}$ is $\sg$
by 2) for $n-1$ and using $\Sigma_1$ uniformization on $(n-1)$ s\.t\. reducts
we can define a good $\sg$ function $\bar f$ s\.t\. $\overline{R}(x,\bar
f(x))\longleftrightarrow\exists\bar y\in H^M_{n-1}\overline{R}(x,\bar y).$
Let $f(x)=h_{n-1}(\bar f(x)).$ Then $f$  is good $\sg$ by 4).

6) \  Let $W$ be universal $\sg$ as in Fact f). By 5) there is a good $\sg$
$g(e,x)$ s\.t\. $\exists y\in\sgn$ hull$(H^M_{n-1}\cup\{x\})$ $W(e,\langle
x,y\rangle)\longleftrightarrow W(e,\langle x,g(e,x)\rangle)$  (and $g(e,x)$
defined $\longrightarrow W(e,\langle x,g(e,x)\rangle)).$  Let
$h_n(e,x)=g(e,x).$  If $y\in\sg$ hull $(\{x\})$ then for some $e, W(e,\langle
x,y'\rangle)\longleftrightarrow y'=y$ so $y=h_n(e,x).$  Clearly
$h_n(e,x)\in\sg$ hull $(\{x\})$ since $h_n$ is good $\sg.$ \hfill{$\dashv$} 

\enddemo

Now we are ready to prove Global $\square.$  Assume $V=L$  and let $\mu$ be a
singular limit ordinal. Our goal is to define $C_\mu,$ a CUB subset of $\mu.$
Let $\beta(\mu)\ge\mu$ be the least limit ordinal $\beta$  such that $\mu$
is not regular with respect to $\widetilde {J}_\beta$-definable functions, and
let $n(\mu)$ be least so that there is a good
$\Sigma^*_{n(\mu)}(\widetilde{J}_{\beta(\mu)})$ {\underbar {partial}}
 function from
an ordinal less than $\mu$ cofinally into $\mu.$  Note that
$\rho^{\beta(\mu)}_{n(\mu)}\le\mu$ as otherwise such a partial function
would belong to $\widetilde{J}_{\beta(\mu)},$ contradicting the leastness of
$\beta(\mu).$  Also $\mu\le\rho^{\beta(\mu)}_{n(\mu)-1},$ else we have
contradicted the leastness of $n(\mu).$ 

For $X\subseteq\widetilde{J}_{\beta(\mu)}$ let $H(X)=\Sigma^*_{n(\mu)}$
hull of $X$  in $\widetilde{J}_{\beta(\mu)}.$  For some least parameter
$q(\mu)\in\widetilde{J}_{\beta(\mu)},$
$H(\mu\cup\{q(\mu)\})=\widetilde{J}_{\beta(\mu)}.$  (``Least'' refers to
the canonical well-ordering of $L.)$  Also let
$\alpha(\mu)=\bigcup\{\alpha<\mu|\alpha=H(\alpha\cup\{q(\mu)\})\cap\mu\}.$
Then (unless $\alpha(\mu)=\bigcup\emptyset=0)$
$\alpha(\mu)=H(\alpha(\mu)\cup\{q(\mu)\})\cap\mu$  and $\alpha(\mu)<\mu.$
To see the latter note that for large enough $\alpha<\mu,
H(\alpha\cup\{q(\mu)\})$ contains both the domain and defining parameter for
a good $\Sigma^*_{n(\mu)}$ partial function from an ordinal less than $\mu$
cofinally into $\mu.$ 

If $\mu<\beta(\mu)$  let $p(\mu)=\langle q(\mu),\mu,\alpha(\mu)\rangle$
and if $\mu=\beta(\mu)$  let $p(\mu)=\alpha(\mu).$ 

We are ready to define $C_\mu.$  Let $C^0_\mu=\{\bar\mu<\mu|$ For some
$\alpha,\bar\mu=\bigcup(H(\alpha\cup\{p(\mu)\})\cap\mu)\}.$ Then $C^0_\mu$
is a closed subset of $\mu.$  If $C^0_\mu$ is unbounded in $\mu$ then set
$C_\mu=C^0_\mu.$  If $C^0_\mu$ is bounded but nonempty then let
$\mu_0=\bigcup C^0_\mu$  and define $C^1_\mu=\{\bar\mu<\mu|$ For some
$\alpha,\bar\mu=\bigcup(H(\alpha\cup\{p(\mu),\mu_0\})\cap\mu)\}.$  If
$C^1_\mu$  is unbounded then set $C_\mu=C^1_\mu.$  If $C^1_\mu$ is bounded
but nonempty then let $\mu_1=\bigcup C^1_\mu$  and define
$C^2_\mu=\{\bar\mu<\mu|$  For some
$\alpha,\bar\mu=\bigcup(H(\alpha\cup\{p(\mu),\mu_0,\mu_1\})\cap\mu)\}.$
Continue in this way, defining $C^k_\mu$  for $k\in\omega$  until $C^k_\mu$
is unbounded or empty. Note that $\alpha_0>\alpha_1>\cdots$  where $\alpha_k$
is greatest so that
$\bigcup(H(\alpha_k\cup\{p(\mu),\mu_0,\cdots,\mu_{k-1}\})\cap\mu )=\mu
_k,$ since $\alpha_k\in
H(\alpha_k\cup\{p(\mu),\mu_0,\cdots,\mu_{k-1},\mu_k\}).$  So for some
least $k(\mu)\in\omega,$  $C^{k(\mu)}_\mu$ is indeed unbounded or empty. If
$C^{k(\mu)}_\mu$  is unbounded then set $C_\mu=C^{k(\mu)}_\mu.$  

If $C^{k(\mu)}_\mu=\emptyset$ then we choose $C_\mu$  to be an
$\omega$-sequence cofinal in $\mu,$  coding approximations to the structure
$\widetilde{J}_{\beta(\mu)},$  as follows. (This is necessary to establish
Relativized $\square$ (c).) Note that
$H=H(\{p(\mu),\mu_0,\cdots,\mu_{k(\mu)-1}\})$ is cofinal in $\mu$  since
$C_\mu^{k(\mu)}=\emptyset.$  Assume first that $n(\mu)=1,$  when $C_\mu$
is more easily described. Then $H$ is also cofinal in $\beta(\mu),$ else
$H\in\widetilde{J}_{\beta(\mu)}$  and $\mu$ is singular inside 
$\widetilde{J}_{\beta(\mu)}.$  Let $h=h_1(e,x)$ be the canonical good
$\Sigma^*_1$  Skolem function for $\widetilde{J}_{\beta(\mu)},$ so 
$H=\{h(e,p)|e\in\omega\}$ where
$p=\{p(\mu),\mu_0,\cdots,\mu_{k(\mu)-1}\}.$  Let
$\bar\sigma_n=\max(\{h(e,p)|e<n\}\cap\mu)$  and
$\sigma_n=\max(\{h(e,p)|e<n\}\cap\beta(\mu)).$  Then
$C_\mu=\{\delta_0,\delta_1,\cdots\}$  where $\delta_n$ is an ordinal coding
TC$(\Sigma^*_1$ hull $(\bar\sigma_n\cup\{p\})$ restricted to $\sigma_n),$
where TC denotes ``transitive collapse''. By the $\Sigma^*_1$ hull of $X$
restricted to $\sigma_n$  we mean the closure of $X$  under $h^{\delta_{n}},$
obtained by interpreting the $\Sigma^*_1$ definition of $h$ in
$\widetilde{J}_{\sigma_{n}}.$  

Now suppose $n(\mu)>1,$ $C^{k(\mu)}_\mu=\emptyset.$  Then if $\rho(\mu)$
denotes $\rho^\mu_{n(\mu)-1},$ $H$ is cofinal in $\rho(\mu),$ else
$H\in\widetilde{J}_{\rho(\mu)}$ and $\mu$ is singular in
$\widetilde{J}_{\rho(\mu)}.$  Let $h$  be the canonical good
$\Sigma^*_{n(\mu)}$  Skolem function for $\widetilde{J}_{\rho(\mu)}$ and let
$p=\{p(\mu),\mu_0,\cdots,\mu_{k(\mu)-1}\}.$  Let
$\bar\sigma_n=\max(\{h(e,p)|e<n\}\cap\mu),$
$\sigma_n=\max(\{h(e,p)|e<n\}\cap\rho(\mu)).$  Then
$C_\mu=\{\delta_0,\delta_1,\cdots\}$  where $\delta_n$ is an ordinal coding
TC$(\Sigma^*_{n(\mu)}$ hull $(\bar\sigma_n\cup\{p\})$ restricted to
$\sigma_n).$  The $\Sigma^*_{n(\mu)}$ hull of $X$  restricted to $\sigma_n$ 
is the closure of $X$  under $h^{\sigma_{n}},$ obtained by replacing the
$(n(\mu)-1)$ st reduct $M_{n(\mu)-1}(x)$  by
$M_{n(\mu)-1}(x)\restriction\sigma_n$  in the $\Sigma^*_{n(\mu)}$ definition
of $h.$  (Recall that $M_n(x)=\langle
J^{A_{n}}_{\omega\rho^{M}_{n}},A_n(x)\rangle;$  by $M_n(x)\restriction\sigma$
we mean $\langle\widetilde{J}^{A_{n}}_\sigma,
A_n(x)\cap\widetilde{J}^{A_{n}}_\sigma\rangle.)$  

Clearly $C_\mu$ is CUB in $\mu,$  and by the same argument used to justify
$\alpha(\mu)<\mu,$  the ordertype of $C_\mu$ is less than $\mu.$  (These
facts are obvious when $C^{k(\mu)}_\mu=\emptyset.)$  So to prove Global
$\square$  we only need to check coherence: $\bar\mu\in\Lim
C_\mu\longrightarrow C_{\bar\mu }=C_\mu\cap\bar\mu.$  

\proclaim{Lemma 4.2}  $\bar\mu\in C^k_\mu\longrightarrow C^k_{\bar\mu
}=C^k_\mu\cap\bar\mu.$  
\endproclaim

\demo{Proof} First suppose that $k=0.$ Given $\bar\mu\in C^0_\mu$ we can
choose $\alpha<\bar\mu$ such that
$\bar\mu=\bigcup(H(\alpha\cup\{p(\mu)\})\cap\mu),$  where $H$  is the
operation of taking the $\Sigma^*_{n(\mu)}$ hull. Also let
$\rho=\bigcup(H(\alpha\cup\{p(\mu)\})\cap\rho^\mu_{n(\mu)-1}).$  Let $\pi:\
\widetilde{J}_{\bar\beta }\longrightarrow\widetilde{J}_{\beta(\mu)}$ be the
inverse to the transitive collapse of $H=\Sigma^*_{n(\mu)}$ hull
$(\bar\mu\cup\{p(\mu)\})$ restricted to $\rho.$   Note that any $x\in H$
belongs to $H(\mu',\rho')=\Sigma^*_{n(\mu )}$ hull $(\mu'\cup\{p(\mu)\})$
restricted to $\rho',$ for some $\mu'<\bar\mu,\rho'<\rho$  and $\mu',\rho'$
can be chosen to be in $H.$  It follows that $H\cap\mu=\bar\mu$ and therefore
when $\mu<\beta(\mu),$ $\pi(\bar\mu)=\mu.$ Also note that 
$\Sigma^*_{n(\mu)-1}$ hull
$(\rho\cup\{p(\mu)\})\cap\rho^\mu_{n(\mu)-1}=\rho,$  so $H$ is closed under
good $\Sigma^*_{n-1}$ functions. It follows that $\pi:\
\widetilde{J}_{\bar\beta }\longrightarrow\widetilde{J}_{\beta(\mu)}$ is
$\Sigma^*_{n-1}$-elementary and $\bar\mu$ is
$\Sigma^*_{n(\mu)-1}(\widetilde{J}_{\bar\beta })$-regular,
$\Sigma^*_{n(\mu)}(\widetilde{J}_{\bar\beta })$-singular. So
$\bar\beta=\beta(\bar\mu),$ $n(\mu)=n(\bar\mu).$  Also
$\pi(q(\bar\mu))=q(\mu).$  Since $\alpha(\bar\mu)<\alpha$ it must be that
$\alpha(\bar\mu)=\alpha(\mu).$  So $\pi(p(\bar\eta))=p(\mu).$  Now it is easy
to see that $C^0_{\bar\mu }=C^0_\mu\cap\bar\mu.$

Now suppose $k=1.$ The above argument shows that  $\bar\mu\in
C^1_\mu\longrightarrow C^0_{\bar\mu }=C^0_\mu\cap\bar\mu$  and hence, since
$\mu_0<\bar\mu,$ $\bar\mu_0=\mu_0.$  Now again, the above argument shows that
$C^1_{\bar\mu }=C^1_\mu\cap\bar\mu.$  The general case $k\ge 0$  now follows
similarly. \hfill{$\dashv$ }

\enddemo

Coherence now follows easily: if $\bar\mu\in\Lim C_\mu$ and $C_\mu=C_\mu^k$
then by Lemma 4.2, $C_{\bar\mu }^k=C_\mu^k\cap\bar\mu$ is unbounded in
$\bar\mu$ so $C_{\bar\mu }=C_{\bar\mu }^k$ and we're done. If
$C_\mu^k=\emptyset$  for some $k$  then lim $C_\mu=\emptyset$ so coherence is
vacuous. 

To establish the appropriate relativized form of $\square$ we need:

\proclaim{Lemma 4.3} Suppose $\pi:\ \langle\widetilde{J}_{\bar\mu },\bar
C\rangle \overset{\Sigma_1}\to\longrightarrow\langle\widetilde{J}_{\mu
},C_\mu\rangle.$  Then $\overline{C}=C_{\bar\mu }$  and $\pi$ extends uniquely
to a $\Sigma^*_{n(\mu)}$-elementary $\tilde\pi:\
\widetilde{J}_{\beta(\bar\mu)}\longrightarrow\widetilde{J}_{\beta(\mu)}$  such
that $p(\mu)\in\Rng\tilde\pi.$ 
\endproclaim

\demo{Proof} First suppose that $C_\mu=C_\mu^k$  for some $k.$  For $\mu'\in
C_\mu$  form $H(\mu')$  as $H$  was formed in the proof of Lemma 2 for
$\bar\mu.$  Then
$\pi(\mu'): \ 
\widetilde{J}_{\beta(\mu')}\longrightarrow\widetilde{J}_{\beta(\mu)}$  
with range $H(\mu')$ is $\Sigma^*_{n(\mu)-1}$-elementary and
$\widetilde{J}_{\beta(\mu)}=\bigcup\{H(\mu')|\mu'\in C_\mu\}.$  And 
$\pi(\mu')\restriction\mu'=id\restriction\mu',$ $\pi(\mu')(p(\mu'))=p(\mu).$
Now let $X=$ Range$(\pi)$  and form $\widetilde{X}=\Sigma^*_{n(\mu)}$ hull
$(X\cup\{p(\mu)\})$ in $\widetilde{J}_{\beta(\mu)}.$  If $y\in\widetilde{X}$
then for some $\mu'\in C_\mu,$  $y=\pi(\mu')(y')$ where
$y'\in\Sigma^*_{n(\mu)}$ hull $((X\cap\widetilde{J}_{\mu'})\cup\{p(\mu')\}).$
In particular if $y\in\widetilde{J}_\mu$ then $y\in\Sigma^*_1$ hull $(X)$ in
$\langle\widetilde{J}_\mu,C_\mu\rangle=X.$  So the inverse to the transitive
collapse of $\widetilde{X}=\tilde\pi$ is a $\Sigma^*_{n(\mu)}$-elementary
embedding extending $\pi,$ with $p(\mu)$ in its range. If $\tilde\pi:\
\widetilde{J}_{\bar\beta }\longrightarrow\widetilde{J}_{\beta(\mu)}$  then
$\bar\mu=\tilde\pi^{-1}(\mu)$ is singular via a
$\Sigma^*_{n(\mu)}(\widetilde{J}_{\bar\beta })$ partial function since either
$\Sigma^*_{n(\mu)}$ hull of $\mu'\cup\{p(\mu)\}$ in
$\widetilde{J}_{\beta(\mu)}$ is unbounded in $\mu$ for some
$\mu'<\bigcup(\Rng(\pi)\cap\mu),$  in which case we can assume
$\mu'\in\Rng\pi$  and by $\Sigma^*_{n(\mu)}$-elementary of $\tilde\pi$ we're
done, or if not $\mu^*=\mu\cap\Sigma^*_{n(\mu)}$ hull $(\mu^*\cup\{p(\mu)\})$
in $\widetilde{J}_{\beta(\mu)}$ where $\mu^*=\bigcup(\Rng\pi\cap\mu),$
contradicting the definition of $\alpha(\mu).$   Since $\bar\mu$ is 
$\Sigma^*_{n(\mu)-1}(\widetilde{J}_{\bar\beta })$-regular, we get 
$\bar\beta=\beta(\bar\mu),$ $n(\mu)=n(\bar\mu).$  Then the
$\Sigma^*_{n(\mu)}$-elementarity of $\tilde\pi$ guarantees that
$\overline{C}=C_{\bar\mu }.$  The uniqueness of $\tilde\pi$  follows from the
fact that $\widetilde{J}_{\beta(\bar\mu)}=\Sigma^*_{n(\mu)}$ hull 
$(\bar\mu\cup\{p(\bar\mu)\})$  and $\tilde\pi\restriction\bar\mu$ is
determined by $\pi.$  

If $C_\mu^k=\emptyset$ for some $k$  then $C_\mu$ was defined as a special
$\omega$-sequence cofinal in $\mu.$  That definition was made precisely to
enable the preceding argument to also apply in this case. \hfill{$\dashv$ }
\enddemo

\flushpar
{\bf Relativized $\bold{\square}$ } \ Let $S=\bigcup\{S_\alpha|\alpha$ an
infinite cardinal$\}.$  There exists  $\langle C_s|s\in S\rangle$ such that
$C_s\in{\Cal{A}}_s$ and:

(a) \  $C_s$ is closed, unbounded in $\mu^0_s,$ ordertype
$(C_s)\le\alpha(s).$ 

If $|s|$ is a successor ordinal then ordertype $(C_s)=\omega.$ 

(b) \ $\nu\in\Lim(C_s)\longrightarrow$ for some $\eta<|s|,$
$\nu=\mu^0_{s\restriction\eta }$ and $C_s\cap\nu=C_{s\restriction\eta }.$

(c) \ Let $\pi:\
\langle\overline{\Cal{A}},\overline{C}\rangle\overset{\Sigma_1}\to
\longrightarrow\langle {\Cal{A}}^0_s,C_s\rangle$ and write
crit$(\pi)=\bar\alpha,$ $\overline{\Cal{A}}=L_{\bar\mu }[\overline{A},\bar
s^*].$ If $\pi(\bar\alpha)=\alpha(s)$ then $L[\overline{A},\bar s^*]\vDash|\bar s|$ is not a cardinal $>\bar
\alpha$  and 

(c1) \ $\overline{C}\in L_\mu [\overline{A},\bar s^*]$  where $\mu$ is the
least p\.r\. closed ordinal greater than $\bar\mu$ s\.t\. $L_\mu
[\overline{A},\bar s^*]\vDash$ card$(|\bar s|)\le\bar\alpha.$

(c2) \ $\pi$ extends to $\pi':\ {\Cal{A}}'\overset{\Sigma_1}\to\longrightarrow
{\Cal{A}}'_s$  where ${\Cal{A}}'=L_{\mu'}[\overline{A},\bar s^*],$ $\mu'=$
largest ordinal either equal to $\bar\mu$ or s\.t\. $\omega\cdot\mu'=\mu'$ and
$L_{\mu'}[\overline{A},\bar s^*]\vDash|\bar s|$ is a cardinal greater than
$\bar\alpha.$

(c3) \ If $\bar\alpha$ is a cardinal and $\pi(\bar\alpha)=\alpha$ then
$\overline{{\Cal{A}}}={\Cal{A}}^0_{\bar s}$  and $\overline{C}=C_{\bar s}.$

\vskip20pt

\flushpar
{\bf Relativized $\bold{\diamond}$ } \  Let $E=$ all $s\in S$ such that
ordertype $(C_s)=\omega.$  There exists $\langle D_s|s\in E\rangle$ such that
$D_s\subseteq {\Cal{A}}^0_s$ and:

(a) \  $D\in\widehat{{\Cal{A}}}_s\ne {\Cal{A}}^0_s,$ $D\subseteq
{\Cal{A}}^0_s\longrightarrow\{\xi<|s|\bigg|s\restriction\xi\in E,
D_{s\restriction\xi }=D\cap{\Cal{A}}^0_{s\restriction\xi }\}$ is stationary in
$\widehat{\Cal{A}}_s.$  

(b) \ $D_s$ is uniformly $\Sigma_1$-definable as an element of ${\Cal{A}}'_s.$

(c) \ If ${\Cal{A}}'_s\vDash\alpha^{++}$ exists then $D_s=\emptyset.$

We now make the necessary modifications to obtain Relativized $\square.$
First, if $\mu$ is a singular limit ordinal and
$\widetilde{J}_\mu\models\alpha$  is the largest cardinal then we thin out
$C_\mu$ to give it ordertype $\le\alpha:$ \ By induction on limit
$\bar\mu\le\mu$ define $C^*_{\bar\mu }$ as follows. For $\bar\mu\le\alpha,$
$C^*_{\bar\mu}=\bar\mu.$  Otherwise $C^*_{\bar\mu }=\{i^{\text{th}}$ element
of 
$C_{\bar\mu }|i\in C^*_{\bar\mu_{0}}$ where $\bar\mu_0=$ ordertype
$(C_{\bar\mu })\}.$  This defines $C^*_{\mu }.$ It is easily verified that
the $C^*_\mu$  enjoy all the properties of the $C_\mu$  except they are only
defined when $\mu$  is a singular limit ordinal such that
$\widetilde{J}_\mu\models$  There is a largest cardinal. In addition,
ordertype $C^*_\mu\le\alpha(\mu),$ the largest cardinal of
$\widetilde{J}_\mu.$  

Now suppose $V=L,$ $\alpha$ is a cardinal, $s\in S_\alpha,$  $|s|>\alpha$ and
$s$  is a $0${\it -string}, meaning that $s(\mu)=0$ for all $\eta\in\Dom(s).$
Then we can choose our predicate $A=\emptyset,$ define $C_s=C^*_{\mu^{0}_{s}}$
and Relativized $\square$ will hold for such $0$-strings. The final comment is
that all we have done will relativize to arbitrary strings $s\in S_\alpha,$
defined relative to an arbitrary predicate $A\subseteq\ORD,\ H_\alpha=L_\alpha
[A]$  for all cardinals $\alpha$

Now we turn to Relativized $\diamond.$  Again we begin with a nonrelativized
version. Let $\alpha$  be a cardinal and assume $V=L.$

\flushpar
{$\underline{\diamond\text{ on }\alpha ^+}$} \  Let $E=$ all $\mu<\alpha^+$
s\.t\. $C_\mu$ has ordertype $\omega.$  There exists $\langle D_\mu|\mu\in
E\rangle$  s\.t\. $D_\mu\subseteq\widetilde{J}_\mu$ and:

(a) \  If $D\subseteq\widetilde{J}_{\alpha^{+}}$ then $\{\mu\in
E|D\cap\widetilde{J}_\mu=D_\mu\}$ is stationary in $\alpha^+.$

(b) \ $D_\mu$ is uniformly $\Sigma_1$ definable as an element of
$\widetilde{J}_{\beta'(\mu)}$  where $\beta'(\mu)=$ largest $\beta$ s\.t\.
either $\beta=\mu$ or $\omega\beta=\beta$ and $\widetilde{J}_\beta\models\mu$
is a cardinal greater than $\alpha.$

(c) \ If $\widetilde{J}_{\beta'(\mu)}\models\alpha^{++}$ exists then
$D_\mu=\emptyset.$  

\demo{Proof} For $\mu\in E$ let $D_\mu=\emptyset$ if
$\widetilde{J}_{\beta'(\mu)}\models\alpha^{++}$ exists and otherwise let 
$\langle D_\mu,F_\mu\rangle$ be least in $\widetilde{J}_{\beta'(\mu)}$ such
that $F_\mu$ is CUB in $\mu$ and $\bar\mu\in F_\mu\longrightarrow\bar\mu\notin
E$ or $D_{\bar\mu }\neq D_\mu\cap\widetilde{J}_{\bar\mu }.$  If
$\langle D_\mu,F_\mu\rangle$ doesn't exist let $D_\mu=\emptyset.$  Properties
(b), (c) are clear. To prove (a), suppose it fails and let $\langle
D,F\rangle$ be least in $\widetilde{J}_{\alpha^{++}}$ such that
$D\subseteq\widetilde{J}_{\alpha^{+}},$  $F$ is CUB in $\alpha^+$ and $\mu\in
F\longrightarrow\mu\notin E$ or $D_\mu\neq D\cap\widetilde{J}_\mu.$
Let $\sigma$  be least such that $\omega\sigma=\sigma$ and $\langle
D,F\rangle\in\widetilde{J}_\sigma.$  Then
$\widetilde{J}_\sigma\models\alpha^+$ is the largest cardinal. Let
$H=\Sigma_1$ Skolem hull of $\{\alpha^+\}$  in $\widetilde{J}_\sigma$  and
$\mu=\bigcup(H\cap\alpha^+).$  Then $\widetilde{J}_{\beta'(\mu)}$ is the
transitive collapse of the $\Sigma_1$ Skolem hull of $\mu\cup\{\alpha^+\}$ in
$\widetilde{J}_\sigma;$ let $\pi:\
\widetilde{J}_{\beta'(\mu)}\longrightarrow\widetilde{J}_\sigma$  have range
equal to the latter hull. Then we have a contradiction provided $\mu\in E.$
But the fact that $H$ is unbounded in $\mu$ implies that $C^0_\mu=\emptyset$
so ordertype $(C_\mu)$ is $\omega$  and $\mu\in E.$ \hfill{$\dashv$ }

\enddemo

Now as with Relativized $\square,$ if $V=L$ and $\alpha$ is a cardinal, $s\in
S_\alpha$  is a $0$-string in $E,$ $|s|>\alpha$  then we can choose
$A=\emptyset$  and define $D_s=D_{\mu^{0}_{s}}$ where the latter comes form
$\diamond$  on $\alpha^+.$  Finally, relativize everything to arbitrary
reshaped strings $s\in S_\alpha$  and an arbitrary predicate $A\subseteq\ORD,$
$L_\alpha [A]=H_\alpha$  for all cardinals $\alpha.$  

This completes the proof of Relativized $\square$  and $\diamond.$ 

\vskip30pt

\centerline{\bf References}

\vskip10pt

\flushpar
Beller-Jensen-Welch \  [82] \, {\it Coding the Universe,} Cambridge University
Press. 

\vskip5pt

\flushpar
Friedman \ [87] \,   Strong  Coding, Annals of Pure and Applied Logic.

\vskip5pt

\flushpar
Friedman \ [94] \,  Jensen's $\Sigma^*$ Theory and the Combinatorial Content
of $V=L,$  to appear, Journal of Symbolic Logic.

\vskip5pt

\flushpar
Jensen \  [72] \,   The Fine Structure of the Constructible Hierarchy, Annals
of Mathematical Logic.

\enddocument